\documentclass{paper}
\usepackage{amssymb,amsmath}
\usepackage{amsfonts}
\input{./defs_allg.tex.input}
\newtheorem{nota}[thm]{Notation}

\begin{document}

\begin{frontmatter}
\title{PBW bases for a class of braided\\ Hopf algebras\thanksref{PhD}}
\thanks[PhD]{This work will be part of the authors PhD thesis written under the supervision of Professor H.-J. Schneider.}
\author{Stefan Ufer\thanksref{BayStaat}}
\address{Mathematisches Institut der Universit\"at M\"unchen, Theresienstr. 39,\\ 80333 M\"unchen, Germany}
\thanks[BayStaat]{Partially supported by Graduiertenf\"orderung des bayerischen Staates and by his parents}
\begin{abstract}
We prove the existence of a basis of Poincar\'e-Birkhoff-Witt type for braided Hopf algebras $R$ generated by a braided subspace $V\subset P(R)$ if the braiding on $V$ fulfils a triangularity condition. We apply our result to pointed Hopf algebras with abelian coradical and to Nichols algebras of low dimensional simple $U_q(\mathfrak{sl}_2)$-modules.
\end{abstract}
\end{frontmatter}

\section{Introduction}
This paper deals with two concepts from combinatorical algebra. The first is the concept of PBW bases.
It is a well known classical result from the theory of Lie algebras that the universal enveloping algebra of a Lie algebra has a basis of elements of the form
\[ x_1^{e_1}\cdot\ldots\cdot x_n^{e_n}\]
with $n\in\N_0,e_1,\ldots,e_n\in\N_0$ and $x_1>\ldots>x_n$ elements of a totally ordered basis of the Lie algebra. This theorem is due to Poincar\'e, Birkhoff and Witt and the basis is called the PBW basis of the universal enveloping algebra. We say that the basis of the Lie algebra generates the PBW basis of the enveloping algebra. We give a formalization of this concept in definition \ref{defi_pbwallg}, incorporating also the analog basis for universal enveloping algebras of $p$-Lie algebras. In \cite{kh.b} Kharchenko proves that such a basis exists for a much larger class of Hopf algebras. Based on his ideas our main theorem proves the existence of a PBW basis for a large class of braided Hopf algebras.
\\[0,2cm]
The other concept is the concept of Lyndon words. These are those words in letters from a set $X$ that are lexicographically smaller than any of their ends (see definition \ref{defi_lyndon}). These words are in bijection with a basis of the free Lie algebra over the set $X$ \cite{Reutenauer}. This bijection offers a close connection between a PBW basis of the enveloping algebra of a Lie algebra $\mathfrak{g}$ generated by some set $X$ (and various relations) and the set of Lyndon words with letters from the set $X$. As the elements generating the PBW basis are iterated Lie brackets in elements from $X$ they can be described as iterated commutators of elements of $X$ in the enveloping algebra \cite{LR}. This correspondence is a very special case of the central theorem in \cite{kh.b} and in this paper.
\\[0,2cm]
More recent results first found by Lusztig \cite{Lusztig_fdha,Lusztig} and Rosso \cite{Rosso_pbw,Rosso_invent}, then also by Ringel \cite{Ringel} showed that also the deformed enveloping algebras of finite dimensional semisimple complex Lie algebras introduced by Drinfel'd \cite{Drin85,Drin87} and Jimbo \cite{Jimbo85} have a PBW basis in the sense of definition \ref{defi_pbwallg}. Later Leclerc observed that the connection to the theory of Lyndon words is preserved under this deformation \cite{Leclerc}. He gives an algorithm to compute the generating set for a PBW basis of the positive part $U^+_q(\mathfrak{g})$ in terms of Lyndon words.
\\[0,2cm]
For general (graded) algebras it is an interesting question whether they admit a (non-trivial) PBW basis or not. The results in the setting of quantum groups treated by now are quite concrete, which is mostly due to the existence of a braid group action \cite{Lusztig}. In \cite{kh.b} Kharchenko follows a more combinatoric approach not relying on a root system or a braid group action. He finds a PBW basis for Hopf algebras generated by an abelian group and a finite set of skew primitive elements such that the adjoint action of the group on the skew primitive generators is given by multiplication with a character. This result can be reinterpreted to provide a PBW basis for braided Hopf algebras in the category of Yetter-Dinfel'd modules over an abelian group which are generated by a completely reducible Yetter-Drinfel'd submodule of primitive elements. In contrast to the results on quantum groups that do only apply to special quotients of the braided tensor algebra, namely the Nichols algebras, Kharchenkos result covers arbitrary braided bialgebra quotients of the tensor algebra.
\\[0,2cm]
It is a priori not clear, if such a strong result generalizes to other situations. In this paper we give a generalization of Kharchenko's result to a bigger family of braided Hopf algebras. The assumption that the braiding on the set of primitive generators is diagonal is replaced by a more general assumption of triangularity (definition \ref{defi_triangular}). The main theorem of this paper says that every braided bialgebra quotient of the tensor algebra of a triangular braided vector space admits a PBW basis. 
\\[0,2cm]
For our combinatorical proof triangularity of the braiding seems to be the natural setting. This situation includes for example a class of braided Hopf algebras constructed from representations of $U_q(\mathfrak{g})$ for $\mathfrak{g}$ a finite dimensional semisimple complex Lie algebra or a Kac-Moody Lie algebra. Remark \ref{rem_pointed} offers a more conceptual characterization of triangular braidings.
\\[0,2cm]
Following a suggestion of the referee we apply this theorem to generalize Kharchenkos existence result for the PBW basis in the sense that the action of the group on the space of skew primitive elements is not required to be via characters any more.
\\[0,2cm]
The approach we take to prove the central theorem is a generalization of Kharchenko's proof and makes intensive use of the comultiplication available in the braided Hopf algebra.
The lack of the braid group action is paid for by a less concrete result. In general the bases given by Kharchenko and in this paper can not be computed without knowledge of the relations in the braided Hopf algebra. However also in the step from diagonal to triangular braidings we lose some information (see remark \ref{rem_losekh}). Nevertheless some results obtained for Nichols algebras of simple $U_q(\mathfrak{sl}_2)$ modules are treated in the last section. They contribute to answering a question raised by N. Andruskiewitsch in \cite{A1}.
\\[0,2cm]

The paper is organized as follows: Section \ref{section_brha} is devoted to the definition and general facts on braided bialgebras. In particular we discuss the Nichols algebra of a braided vector space. Section \ref{section_lyndon} provides basic facts about Lyndon words. For proofs of the results in this section the reader is referred to the literature \cite{Lothaire,Reutenauer,Ufnarovski}. In Section \ref{section_brcomm} the setting for the central theorem - braidings of left (resp. right) triangular type - are introduced. We study braided commutators in the tensor algebra of the braided vector space. These will be a major tool for the proof of the main theorem. Section \ref{section_comult} deals with technical combinatorical facts about the comultiplication in the tensor algebra of a left triangular braided vector space. These results rely heavily on the condition of triangularity. In section \ref{section_PBWbasis} we define the data describing the PBW basis for any braided bialgebra quotient of the tensor algebra and prove the main theorem. Furthermore we prove a result providing an important restriction on the heights of the PBW generators. Section \ref{sect_right} offers a transfer of the results to right (instead of left) triangular braidings. Section \ref{sect_applHA} contains the application of our result to Hopf algebras generated by an abelian group and finitely many skew primitive elements. Finally in section \ref{section_examples} we study Nichols algebras of simple \Uq{sl_2} modules of low dimension.\par
Throughout the paper $k$ will be a field, all vector spaces will be $k$-vector spaces and all tensor products are taken over $k$.
\\[0,2cm]
I would like to thank the referee for useful comments, my supervisor Prof. Schneider for his guidance and Istv\'an Heckenberger for helpful discussions.
\newpage

\section{Braided vector spaces and braided bialgebras}
\label{section_brha}
In this section we review some basic facts about braided Hopf algebras.

\begin{defn}
A braided vector space is a vector space $V$ together with an automorphism $c$ of $V\otimes V$ that satisfies the braid equation:
\[ (c\otimes \id_V)(\id_V\otimes c)(c\otimes \id_V) = (\id_V\otimes c)(c\otimes \id_V)(c\otimes \id_V).\]
\end{defn}

We define further isomorphisms
\[ c_{n,m}:V^{\otimes n}\otimes V^{\otimes m} \rightarrow V^{\otimes m}\otimes V^{\otimes n} \mtxt{for} m,n\geq 0\]
inductively by $c_{0,1}=\id_V,c_{1,0}=\id_V,c_{0,0}=\id_k,c_{1,1}=c$ and
\[c_{1,m+1} = (\id_V\otimes c_{1,m})(c\otimes \id_{V^{\otimes m}}), c_{n+1,m}= (c_{n,m}\otimes \id_V)(\id_{V^{\otimes n}}\otimes c_{1,m}). \]
This induces an automorphism $\hat{c}:T(V)\otimes T(V) \rightarrow T(V)\otimes T(V)$ turning the tensor algebra $T(V)$ of $V$ into a braided vector space $(T(V),\hat{c})$. From now on we will denote the braiding on $T(V)$ by $c$ as well.\newline
If $(V,c)$ is a braided vector space and $f:V^{\otimes n}\rightarrow V^{\otimes m}$ is a homomorphism we say that $f$ \textsl{commutes} with $c$ if
\[ (f\otimes \id_V)c_{1,n} = c_{1,m}(\id_V\otimes f) \mtxt{and} (\id_V\otimes f)c_{n,1} = c_{m,1}(f\otimes \id_V).\]
The braid equation says that $c$ commutes with $c$.\newline
Assume $A$ is an algebra and $(A,c)$ a braided vector space such that the product $m:A\otimes A\rightarrow A$ and the unit $\eta:A\rightarrow k$ of $A$ commute with $c$. Define $A\obar A:=A\otimes A$ as a vector space. Endowed with the maps
\begin{eqnarray*}
&m_{A\obar A}& := (m_A\otimes m_A)(\id_A\otimes c\otimes \id_A): (A\obar A)\otimes (A\obar A)\rightarrow A\obar A,\\
&\eta_{A\obar A}&:= \eta_A\otimes\eta_A :k\rightarrow A\obar A.
\end{eqnarray*}
this is an algebra. Furthermore the opposite algebra $A^{op,c}$ is the algebra obtained from $A$ using the multiplication $m_{A^{op,c}}=m_A c$.\newline
Dually, if $C$ is a coalgebra and $(C,c)$ is a braided vector space such that $\Delta,\eps$ commute with $c$, define $C\obar C:=C\otimes C$ as a vector space. Endowed with the maps
\begin{eqnarray*}
&\Delta_{C\obar C}&:= (\id_C\otimes c\otimes \id_C)(\Delta_C\otimes\Delta_C):C\obar C\rightarrow (C\obar C)\otimes(C\obar C),\\
&\eps_{C\obar C}& := \eps_C\otimes\eps_C: C\obar C\rightarrow k.
\end{eqnarray*}
this is a coalgebra. The following definition of a braided Hopf algebra is equivalent to that of \cite{Takeuchi_survey}.

\begin{defn}
A braided bialgebra $(R,c)$ (or $R$) is a tuple $(R,m,\eta,\Delta,\eps,c)$ such that
\begin{itemize}
\item $(R,m,\eta)$ is an algebra,
\item $(R,\Delta,\eps)$ is a coalgebra,
\item $(R,c)$ is a braided vector space,
\item $m,\eta,\Delta,\eps$ commute with $c$
\end{itemize}
and one of the two following equivalent conditions holds:
\begin{itemize}
\item $\Delta:R\rightarrow R\obar R$ and $\eps:R\rightarrow k$ are algebra homomorphisms.
\item $m:R\obar R\rightarrow R$ and $\eta:k\rightarrow R$ are coalgebra homomorphisms.
\end{itemize}
A homomorphism $f:(R,c)\rightarrow (R',d)$ of braided bialgebras is a homomorphism of algebras and coalgebras such that $(f\otimes f)c = d(f\otimes f)$.\newline
As usual $\Hom_k(R,R)$ is an algebra with the convolution product and $R$ is called a braided Hopf algebra if the identity on $R$ is convolution invertible. In this case the convolution inverse of the identity is called the antipode of $R$.\newline
An element $x\in R$ is called \emph{primitive} if we have $\Delta(x) = 1\otimes x + x\otimes 1$. Let $P(R):= \{x\in R|x\;\mbox{primitive}\}$.
\end{defn}

If $(R,c)$ is a braided Hopf algebra, the antipode commutes with the braiding \cite{Takeuchi_survey}.
\begin{rem}
Let $H$ be a Hopf algebra with bijective antipode.
Every bialgebra (Hopf algebra) in the category of Yetter-Drinfel'd modules over $H$ is a braided bialgebra (braided Hopf algebra) in the sense of this definition. Conversely Takeuchi shows that every rigid braided bialgebra (Hopf algebra) can be realized as a bialgebra (Hopf algebra) in the category of Yetter-Drinfel'd modules over some Hopf algebra with bijective antipode \cite{Takeuchi_survey}.
\par
Nevertheless our notion of a morphism of braided bialgebras is weaker than that of a morphism of bialgebras in a Yetter-Drinfel'd category. Assume that we have a bialgebra $R$ in the category of Yetter-Drinfel'd modules over some Hopf algebra. A subbialgebra $R'$ in this setting is a Yetter-Drinfel'd submodule and thus we have automatically induced braidings
\[ R'\otimes R\rightarrow R\otimes R', R'\otimes R\rightarrow R\otimes R' \mtxt{and} R'\otimes R'\rightarrow R'\otimes R'.\]
On the other hand assume we have a braided bialgebra $R''$ that is a braided subbialgebra of $R$ in the sense that the inclusion is a morphism of braided bialgebras, but $R''$ is not necessarily a Yetter-Drinfel'd submodule. In this case we obtain only a braiding for $R''$
\[R''\otimes R'' \rightarrow R''\otimes R''.\]
In \cite{Takeuchi_survey} $R''$ is called a \textsl{non-categorical} (braided) subbialgebra of $R$ in this case.
\end{rem}
\begin{lem}
Let $(V,c)$ be a braided vector space. Then the tensor algebra $T(V)$ admits a unique structure of a braided Hopf algebra such that the elements of $V$ are primitive.
\end{lem}
 \begin{pf}
Define $\Delta$, $\eps$ and $S:T(V)\rightarrow T(V)^{op,c}$ using the universial property of the tensor algebra such that for all $v\in V$ we have $\Delta(v)=1\obar v+v\obar 1,\eps(v)=0,S(v)=-v$. So $\Delta$ and $\eps$ are algebra homomorphisms. It is easy to check that the set of elements $x\in T(V)$ satisfying $(\id\star S)(x)=\eps(x)$ and $(S\star \id)(x) = \eps(x)$ is closed under multiplication. Checking this equation for the generators is trivial, so $S$ is a convolution inverse of $\id$. Of course $\eps$ commutes with $c$. By the construction, $\Delta$ is a composition of homomorphisms of the type $V^{\otimes i}\otimes c\otimes V^{\otimes j}, i,j\geq 0$ (for exact formulas see \cite{Sbg_borel}). As $c$ commutes with $c$, these homomorphisms commute with $c$ and so does $\Delta$.
 \qed\end{pf}

\begin{rem}
If $(R,c_R)$ is a braided bialgebra and $V:=P(R)$ its space of all primitive elements, we have
\[ c_R(V\otimes V)\subset V\otimes V.\]
\end{rem}

\begin{lem}
\label{epiforbraidedHA}
Let $(R,c_R)$ be a braided bialgebra, $V\subset P(R)$ a braided subspace (this means $c_R(V\otimes V) = V\otimes V$). Then there is a unique homomorphism of braided bialgebras $\pi:T(V,c_R|V\otimes V)\rightarrow R$ with $\pi|V=\id_V$.
\end{lem}
 \begin{pf}
Uniqueness is obvious. Of course $\pi$ exists as algebra homomorphism. Denote the braiding on the tensor algebra induced by $c_R|V\otimes V$ by $c_{T(V)}$. Using the universial property of the tensor algebra we obtain that $\pi$ is a coalgebra homomorphism, provided $\pi\otimes\pi:T(V)\obar T(V)\rightarrow R\obar R$ is an algebra homomorphism. It is easy to check this, if $(\pi\otimes\pi)c_{T(V)} = c_R(\pi\otimes\pi)$. So we are left to show this. By construction we have $\pi|V^{\otimes l} = m_l|V^{\otimes l}$, a restriction of the $l$-fold multiplication of $R$. Thus for all $r,s\geq 0$
 \begin{eqnarray*}
 (\pi\otimes\pi)c_{T(V)}|V^{\otimes r}\otimes V^{\otimes s} &=& (m_s\otimes m_r){(c_{T(V)})}_{r,s}|V^{\otimes r}\otimes V^{\otimes s}\\ 
 &=& c_R (m_r \otimes m_s)|V^{\otimes r}\otimes V^{\otimes s} \\
 &=& c_R(\pi\otimes\pi)|V^{\otimes r}\otimes V^{\otimes s},
 \end{eqnarray*}
 where the second equality is because the multiplication of $R$ commutes with $c$.
 \qed\end{pf}

\begin{defn}
Let $(R,c)$ be a braided bialgebra. A subspace $I\subset R$ is called a (braided) biideal, if it is an ideal, a coideal and 
\[ c(R\otimes I + I\otimes R) = R\otimes I + I\otimes R.\]
If $(R,c)$ is a braided Hopf algebra with antipode $S$, $I$ is called a (braided) Hopf ideal if it is a biideal with $S(I)\subset I$.
\end{defn}
\begin{lem}
Let $(R,c)$ be a braided bialgebra.
\begin{enumerate}
\item If $I\subset R$ is a braided biideal there is a unique structure of a braided bialgebra on the quotient $R/I$ such that the canonical map is a homomorphism of braided bialgebras.
\item If $\pi:(R,c)\rightarrow(S,d)$ is a morphism of braided bialgebras, $\ker\pi$ is a braided biideal of $R$.
\item Analogous statements hold for braided Hopf ideals.
\end{enumerate}
\end{lem}
 \begin{pf}
 Part 1: Uniqueness is clear because $\pi$ is surjective. Obviously $R/I$ is an algebra and a coalgebra in the usual way with structure maps $\bar{m},\bar{\eta},\bar{\Delta}$ and $\bar{\eps}$. Furthermore $c(\ker(\pi\otimes\pi))=\ker(\pi\otimes\pi)$ and thus $c$ induces an automorphism $\bar{c}$ of $R/I\otimes R/I$ such that $(\pi\otimes\pi)c = \bar{c}(\pi\otimes\pi)$. Surjectivity of $\pi\otimes\pi$ ensures that $\bar{c}$ satisfies the braid equation and that $\bar{\Delta},\bar{\eps}$ are algebra homomorphisms. $\bar{m},\bar{\eta},\bar{\Delta},\bar{\eps}$ commute with $\bar{c}$ because $m,\eta,\Delta,\eps$ commute with $c$ and $\pi$ is surjective.\newline
 Part 2: Of course $I:= \ker\pi$ is an ideal and a coideal. It remains to show that the condition for $c$ holds. As $(\pi\otimes\pi)c=d(\pi\otimes\pi)$ and $c$ is bijective we have $c(\ker(\pi\otimes\pi))=\ker(\pi\otimes\pi)$. In view of $\ker(\pi\otimes\pi)= I\otimes R+R\otimes I$ the proof is complete.
 \qed\end{pf}

\begin{exmp}
\label{defi_nichols}
\emph{The Nichols algebra of a braided vector space.} Let $(V,c)$ be a braided vector space. The Nichols algebra $\mathcal{B}(V,c)$ as defined in \cite{AS5} satisfies
\begin{itemize}
\item $\mathcal{B}(V,c)=\bigoplus_{n\in\N_0}R(n)$ is a graded braided Hopf algebra (this means graded as algebra, coalgebra and braided vectorspace simultaneously),
\item $R(0)\isom k$ and $R(1)\isom V$ as braided vector spaces,
\item $R(1) = P(R)$ and
\item $R$ is generated by $R(1)$ as an algebra.
\end{itemize}
The Nichols algebra for a braided vector space $(V,c)$ exists and is unique up to isomorphism. It is a braided analogue of the classical symmetric algebra.
See \cite{Sbg_borel,AS5} for more details.
\end{exmp}

\section{Lyndon words}
\label{section_lyndon}
The PBW basis constructed later is closely related to Lyndon words in letters from a set $X$ of primitive elements of the braided bialgebra. In this section we discuss the necessary facts about Lyndon words.
Let ($X$,$\lx$) be a finite totally ordered set and $\X$ the set of all words in the letters $X$ (the free monoid over $X$). Throughout this section all words will be from $\X$. Recall that the \textsl{lexicographical order} on $\X$ is the total order defined in the following way: For words $u,v\in\X,\, u\lx v$ iff either $v\in u\X$ ($u$ is the beginning of $v$) or if there exist $r,s,t\in\X, a,b\in X$ such that
\[ u=ras, v=rbt \mtxt{and} a<b.\]
For example if $x,y\in X,x\lx y$ then $x\lx xy\lx y$.

\begin{nota}
For a word $u\in \X$ let $l(u)$ be the length of $u$. Define $\X^n := \{u\in\X | l(u)=n\}$, for $v\in\X$ let $\X_{\gx v} := \{u\in\X |u\gx v\}, \X_{\gex v} := \{u\in\X |u\gex v\}$, $\X_{\gx v}^n := \X^n\cap\X_{\gx v}$ and $\X_{\gex v}^n:=\X^n\cap\X_{\gex v}$.
\end{nota}

\begin{defn}
\label{defi_lyndon}
Let $u\in\X$. The word $u$ is called a Lyndon word if $u\neq 1$ and $u$ is smaller than any of its proper endings. This means for all $v,w\in\X\setminus\{1\}$ such that $u=vw$ we have $u\lx w$.\newline
These words are also called regular words in \cite{Ufnarovski} or standard words in \cite{kh.b}.
\end{defn}
A word $u$ is Lyndon if and only if for every factorization $u=vw$ of $u$ into nonempty words $v,w$ we have $u=vw\lx wv$ (\cite{Lothaire}, 5.1.2.)
\begin{thm}(Lyndon, see \cite{Lothaire}, theorem 5.1.5.)\newline
Any word $u\in\X$ may be written uniquely as a nonincreasing product of Lyndon words
\[u=l_1l_2\ldots l_r, \mtxt{$l_i$ Lyndon words and} l_1\gex l_2\gex\ldots\gex l_r.\]
\end{thm}
This decomposition is obtained inductively by choosing $l_1$ to be the longest beginning of $u$ that is a Lyndon word.
It will be refered to as the \textsl{Lyndon decomposition} of $u$. The occuring Lyndon words are called the \textsl{Lyndon letters} of $u$.
\begin{thm}(\cite{Reutenauer}, theorem 5.1. and section 4.1.)\newline
The set of Lyndon words is a Hall set with respect to the lexicographical order. This means that for every Lyndon word $u\in\X\setminus X$ we have a fixed decomposition $u=u'u''$ into nonempty Lyndon words $u',u''$ such that either $u'\in X$ or the decomposition of $u'$ has the form $u'=vw$ with $w\gex u''$.
\end{thm}
This decomposition is obtained by choosing $u''$ to be the minimal (with respect to the lexicographical order) or (equivalently) the longest end of $u$ that is Lyndon. As in \cite{kh.b} it is referred to as the \textsl{Shirshov decomposition} of $u$.
\begin{lem}(\cite{kh.b}, Lemma 5)\newline
For $u,v\in\X$ we have $u\lx v$ if and only if $u$ is smaller than $v$ when comparing them using the lexicographical order on the Lyndon letters. This means if $v=l_1\ldots l_r$ is the Lyndon decomposition of $v$, we have $u\lx v$ iff
\begin{itemize}
\item $u$ has Lyndon decomposition $u=l_1\ldots l_i$ for some $0\leq i < r$
\item or $u$ has Lyndon decomposition $u=l_1\ldots l_{i-1}\cdot l\cdot l'_{i+1}\ldots l'_s$ for some $1\leq i < r,s\in\N$ and some Lyndon words $l,l'_{i+1},\ldots,l'_s$ with $l\lx l_{i}$.
\end{itemize}
\end{lem}

\section{Braided commutators in the tensor algebra}
\label{section_brcomm}
Major tools for constructing the PBW Basis will be the braided commutators discussed in this section.
Take a finite dimensional vectorspace $V$, an endomorphism $r$ of $V\otimes V$ satisfying the braid equation and a basis $X$ of $V$. Define the endomorphism $r_{n,m}:V^{\otimes n}\otimes V^{\otimes m} \rightarrow V^{\otimes m}\otimes V^{\otimes n}$ in the same way as for braidings. We will omit the indices $n,m$ whenever it is clear from the context which endomorphism is used. \par
Identify $k\X$ - the free algebra over $X$ - with the tensor algebra of $V$ in the obvious way and construct a $k$-linear endomorphism $[-]_r$ of $k\X$ inductively. Set for all $x\in X$ 
\[[1]_r:=1 \mtxt{and} [x]_r:=x.\]
For Lyndon words $u\in\X$ of degree $>1$ with Shirshov decomposition $u=vw$ define 
\[[u]_r:=m(\id-r_{l(v),l(w)})([v]_r\otimes [w]_r),\]
 where $m$ denotes multiplication in $k\X$. For arbitrary words with Lyndon decompositon $u=u_1\ldots u_t$ let 
\[[u]_r := [u_1]_r\ldots[u_t]_r.\]
Obviously $[-]_r$ is a graded homomorphism of the graded vectorspace $k\X$. The idea of using a homomorphism of this type to construct PBW bases can be found in \cite{kh.b} and is motivated by the theory of (free) Lie algebras.

For the rest of this paper we will only deal with braidings that fulfil a certain combinatorical condition.
\begin{defn}
\label{defi_triangular}
Let $V$ be a finite dimensional vector space with a totally ordered basis $X$ and $c\in\End(V\otimes V)$.\newline
The endomorphism $c$ will be called \emph{left triangular} (with respect to the basis $X$) if for all $x,y,z\in k$ with $z\gx y$ there exist $\gamma_{x,y}\in k$ and $v_{x,y,z}\in V$ such that for all $x,y\in X$
\[ c(x\otimes y) = \gamma_{x,y} y\otimes x + \sum\limits_{z\gx y} z\otimes v_{x,y,z} . \]
The endomorphism $c$ will be called \emph{right triangular} (with respect to the basis $X$) if for all $x,y,z\in k$ with $z\gx y$ there exist $\beta_{x,y}\in k$ and $w_{x,y,z}\in V$ such that for all $x,y\in X$
\[ c(x\otimes y) = \beta_{x,y} y\otimes x + \sum\limits_{z\gx x} w_{x,y,z}\otimes z . \]
A braided vector space $(V,c)$ will be called left (resp. right) triangular with respect to the basis $X$ if $c$ is left (resp. right) triangular with respect to the basis $X$.
\end{defn}

\begin{rem}
The name "left triangular" is motivated by the following observation: Assume in the situation of the definition that $V$ has dimension $n$ and denote by $B=(b_1,\ldots,b_{n^2})$ the basis $\{x\otimes y|x,y\in X\}$ of $V\otimes V$ ordered lexicographically. By $B^{op} = (b'_1,\ldots,b'_{n^2})$ denote the basis obtained from $B$ by flipping the sides of every tensor (not changing the order). Then the matrix $A\in\GL(n^2,k)$ satisfying $c(b'_1,\ldots,b'_{n^2}) = (b_1,\ldots,b_{n^2})A$ has the following form:
\[ A = \Matrix{cccc}{D_1&0&\hdots&0\\\star&D_2&\ddots&\vdots\\ \vdots&\ddots&\ddots&0\\\star&\hdots&\star&D_n},\]
where $D_1,\ldots,D_n\in\GL(n,k)$ are diagonal matrices. If the braiding was diagonal (defined after the next example), this matrix would be diagonal.
\end{rem}

\begin{rem}
\label{rem_pointed}
In \cite{Ich_pointed} we show that for a braided vectorspace $(V,c)$ the following statements are equivalent:
\begin{itemize}
\item $c$ is right triangular
\item There is a pointed Hopf algebra $H$ with abelian coradical having $M$ as a (left-left) Yetter-Drinfel'd module such that the induced braiding is $c$ and $G(H)$ acts diagonally on $M$.
\end{itemize}
For technical reasons we will first deal with left triangular braidings. In section \ref{sect_right} we will transfer our results to the right triangular case.
\end{rem}

\begin{exmp}
\label{expl_uqg}
Assume $k=\C$.
Let $\mathfrak{g}$ be a finite dimensional semisimple complex Lie algebra and $q\in\C$ not a root of unity. Let $\Phi\subset V$ be the root system of $\mathfrak{g}$ ($V$ a finite dimensional $\Q$ vector space), fix a basis $\Pi = \{\mu_1,\ldots,\mu_s\}$ of it and denote by $\Lambda\subset V$ the weight lattice. For any finite dimensional \Uq{g}-module $M$ of type $1$ with weight space decomposition $M=\bigoplus_{\mu\in\Lambda}M_\mu$ we have a \Uq{g} linear braiding (see \cite{Jantzen}, Chapter 7)
\[\Theta:M\otimes M\rightarrow M\otimes M.\]
This braiding is built up of $\C$-linear maps
\[ \Theta_\nu:M\otimes M\rightarrow M\otimes M \mtxt{for} \nu\in\Z\Phi,\nu\geq 0\]
satisfying for all $\mu,\mu'\in\Lambda,\nu\in\Z\Phi,\nu\geq 0$
\[\Theta_\nu(M_\mu\otimes M_{\mu'})\subset M_{\mu'-\nu}\otimes M_{\mu+\nu}\]
and a map
\[ f:\Lambda\times\Lambda\rightarrow\C\]
in the following way: For $a\in M_\mu,b\in M_{\mu'}$ we have
\[ \Theta(a\otimes b) = f(\mu',\mu)\lsum{\nu\geq 0}{} \Theta_\nu(b\otimes a).\]
Note that the sum is actually finite because for finite dimensional modules $M$ only finitely many weight spaces $M_\mu$ can be different from zero.
We will construct a basis $B$ of $M$ such that the braiding $\Theta$ is left triangular with respect to this basis.\par
Now consider the total order $\gtr$ defined on $V$ using the basis $\Pi$ in the following way:
\[ \lsum{i=1}{s}a_i\mu_i\gtr\lsum{i=1}{s}b_i\mu_i \Leftrightarrow (a_1,\ldots,a_s)>(b_1,\ldots,b_s),\]
where on the right side we order the sequences in $\Q^s$ lexicographically by identifying them with words of $s$ letters from $\Q$.
\par
Then for $\mu,\mu',\nu\in V$ $\mu\gtr\mu'$ implies $\mu+\nu\gtr\mu'+\nu$ and $\nu>0$ implies $\nu\gtr 0$. For every $\mu\in\Lambda$ with $M_\mu\neq 0$ choose a totally ordered basis $(B_\mu,\lex)$ of $M_\mu$ and order the union $B = \cup_\mu B_\mu$ by requiring that for $b\in B_\mu,b'\in B_{\mu'}$
\[ b\lx b' \Leftrightarrow \mu\gtr\mu'.\]
This defines a totally ordered basis of $M$ and for $b\in B_\mu,b'\in B_{\mu'}$ we have (using $\Theta_0 = id$)
\[\Theta(b\otimes b') = f(\mu',\mu) \left(b'\otimes b + \lsum{\nu>0}{}\Theta_v(b'\otimes b)\right)
\in     f(\mu',\mu) b'\otimes b + \lsum{\nu\ltr\mu'}{} M_\nu\otimes M\]
showing that the braiding is indeed left triangular. In the same one sees that the braiding is also right triangular.
\end{exmp}
We thank the referee for pointing out the following example.
\begin{exmp}
\label{exmp_jordan}
Assume that $k$ is algebraically closed.
Let $G$ be an abelian group and $V\in{}^G_G\mYD$ a Yetter-Drinfel'd module over $G$. Then the induced braiding
\[ c:V\otimes V\rightarrow V\otimes V,\:\: c(v\otimes w) = v\sm1 w\otimes v\s0 \]
is left triangular.
\end{exmp}
\begin{pf}
For all $g\in G$ let $V_g := \{v\in V|\delta(v) = g\otimes v\}$. Then the $V_g$ are $G$-submodules of $V$. As every simple submodule of a finite dimensional $G$ module is one dimensional we see that each $V_g$ has a flag of invariant subspaces. So for all $g\in G$ we find a basis $v^g_1,\ldots,v^g_{r_g}$ of $V_g$ such that for all $h\in G$
\[h\cdot v^g_i\in kv^g_i\oplus\ldots\oplus kv^g_{r_g}.\]
Now by concatenating these bases and ordering each according to the indices we obtain a totally ordered basis such that $c$ is triangular.
\qed\end{pf}

In the previous example we used
\begin{lem}
Let $G$ be an abelian group and $V$ a finite dimensional simple representation of $G$. Then $V$ is onedimensional.
\end{lem}
\begin{pf}
Let $V$ be a finite dimensional simple $G$ module. By Schurs lemma (and because $k$ is algebraically closed) there are scalars $\lambda_g\in k\setminus\{0\}$ such that for all $g\in G,v\in V$
\[ gv = \lambda_g v.\]
So every $v\in V$ spans a onedimensional submodule. As $V$ is simple we have that $V$ is onedimensional.
\qed\end{pf}

\par
Trivial but useful examples of left (and right) triangular endomorphisms are \emph{diagonal braidings}. These are braidings $d$ of $V\otimes V$ such that there is a basis $X\subset V$ and coefficients $\alpha_{x,y}\in k$ satisfying for all $x,y\in X$.
\[ d(x\otimes y) = \alpha_{x,y} y\otimes x \]
In this case we have for arbitrary words $u,v\in\X$
\[ d(u\otimes v) = \alpha_{u,v} v\otimes u,\]
where the coefficients $\alpha_{u,v}$ are defined inductively for all $x\in X$ by
\[ \alpha_{1,1} = 1\, , \, \alpha_{x,1}=1=\alpha_{1,x}\]
and for all $u,u',v,v'\in \X$ by
\[ \alpha_{uu',v} = \alpha_{u,v}\alpha_{u',v} \mtxt{resp.} \alpha_{u,vv'}=\alpha_{u,v}\alpha_{u,v'}.\]
Assume $V$ is a vector space and $r$ an endomorphism of $V\otimes V$ that is left triangular with respect to a basis $X\subset V$. Define the endomorphism $d:V\otimes V\rightarrow V\otimes V$ by 
\[d(x\otimes y) = \gamma_{x,y} y\otimes x \mtxt{for all} x,y\in X,\]
where the coefficients $\gamma_{x,y}$ are those occuring in definition \ref{defi_triangular} for $r$. This endomorphism is a braiding and is called the \textsl{diagonal component} of $r$.

\begin{rem}
\label{remark_notdiag}
There are braidings that are triangular but not diagonal. For example the braiding on the simple twodimensional \Uq{sl_2} module $(M,c)$ is left and right triangular, but not diagonal. Observe that if $c$ was diagonal with respect to some basis $A$ and diagonal coefficients $\alpha_{a,b},a,b\in A$, then $c$ would be diagonalizeable as endomorphism of $M\otimes M$ with eigenvalues $\pm \sqrt{\alpha_{a,b}\alpha_{b,a}}$ for $b\neq a$ (eigenvectors $\sqrt{\alpha_{b,a}}a\otimes b \pm \sqrt{\alpha_{a,b}}b\otimes a$) resp. $\alpha_{a,a}$. But the eigenvalues of $c$ in our case are $-1$ and $q^{-2}$. As we assumed $q$ not to be a root of unity, the braiding can not be diagonal.
\end{rem}

In \cite{kh.b} the case of diagonal braidings is studied. The central problem of this section is to generalize results of the diagonal case and to provide new tools necessary in the triangular case. The next lemma for example is trivial in the diagonal case.
\begin{lem}
Let $V$ be a vector space and assume that $c$ is a left triangular endomorphism with respect to the basis $X$. We have for words $u,v\in \X$:
\[c(u\otimes v) \in d(u\otimes v) + k\X^{l(v)}_{\gx v}\otimes k\X^{l(u)},\]
where $d$ is the diagonal component of $c$.
\end{lem}
\begin{pf}
We use double induction on $l(u)$ and $l(v)$. For $l(u)=0$ and for $l(v)=0$ the claim is trivial. So from now on assume $l(u),l(v)>0$. If $l(u)=l(v)=1$ the claim is exactly the condition from definition \ref{defi_triangular}. Now let $l(u)=1,l(v)>1$ and write $v=xw$ with $x\in X,w\in\X$. Use the notation from definition \ref{defi_triangular}. Then with $q:=l(v)$ the induction hypothesis gives
\begin{eqnarray*}
c_{1,q}(u\otimes v) &=& (\id_V\otimes c_{1,q-1})(c_{1,1}(u\otimes x)\otimes w)\\
&=& \gamma_{u,x} (\id\otimes c_{1,q-1})((x\otimes u)\otimes w) \\
&&+\lsum{z\gx x}{} (\id\otimes c_{1,q-1})((z\otimes v_{u,x,z})\otimes w)\\
&\in& \gamma_{u,x}\gamma_{u,w} x\otimes w\otimes u + \lsum{z\gx x}{} z\otimes k\X^{q-1}\otimes k\X^1\\
&\subset& \gamma_{u,xw} xw\otimes u + k\X^q_{\gx xw}\otimes k\X^1,
\end{eqnarray*}
where the last inclusion follows from the definition of the lexicographical order (note that in any case only words of the same length are compared). So now assume $q=l(v)\geq 1,p:=l(u)>1$ and write $u=wx$ for some $x\in X$. Then
\begin{eqnarray*}
c_{p,q}(u\otimes v) &=& (c_{p-1,q}\otimes \id_V)(w\otimes c_{1,q}(x\otimes v))\\
&\in& \gamma_{x,v} c_{p-1,q}(w\otimes v)\otimes x + c_{p-1,q}(w\otimes k\X^q_{\gx v}) \otimes k\X^1\\
&\subset& \gamma_{x,v}\gamma_{w,v} v\otimes wx + k\X^q_{\gx v}\otimes k\X^{p-1}\otimes x + k\X^q_{\gx v}\otimes k\X^p
\end{eqnarray*}
using the induction hypothesis for $p$ twice.
\qed\end{pf}

\begin{nota}
Let $(V,c)$ be a braided vector space that is left triangular with respect to a basis $X$. An endomorphism $r$ of $V\otimes V$ will be called \textsl{admissible} if it satisfies the braid equation and is left triangular with respect to the basis $X$.
\end{nota}
For example the braiding $c$ itself, braidings which are diagonal with respect to the basis $X$ and the zero morphism are admissible. The concept of commutators induced by admissible endomorphism allows us to formulate the process Kharchenko \cite{kh.c} refers to as monomial crystallization, namely the transfer from a basis of iterated commutators to a basis made up of the underlying words. The first part of the following lemma is a generalization of the second part of \cite[Lemma 5]{kh.b} to our case of commutators comming from arbitrary admissible endomorphisms.

\begin{lem}
\label{lemma_smallestterm}
Let $(V,c)$ be a left triangular braided vector space with basis $X$ and $r$ an admissible endomorphism. Then for every word $u\in\X$ the polynomial $[u]_r$ is homogenous of degree $l(u)$ and the smallest monomial in this term is $u$ with coefficient $1$:
\[ [u]_r\in u + k\X^{l(u)}_{\gx u}. \]
In particular if the diagonal component of the braiding $c$ has the coefficients $\gamma_{x,y}$ and $r$ is itself diagonal, we have
\[ c([u]_r\otimes [v]_r)\in \gamma_{u,v} [v]_r\otimes [u]_r + k\X^{l(v)}_{\gx v}\otimes k\X^{l(u)}.\]
\end{lem}
\begin{pf}
Proceed by induction on $l(u)$. The cases $l(u)=0,1$ follow from the definition of $[-]_r$. In the case $l(u)>1$ first assume $u$ is a Lyndon word. Then we have a Shirshov decomposition $u=vw$ of $u$. With $p:=l(v),q:=l(w)$ ($m$ is the multiplication map) we have
\[ [u]_r = [v]_r[w]_r - m r_{p,q}([v]_r\otimes [w]_r)\]
and using the induction assumption we obtain
\begin{eqnarray*}
[u]_r &\in& (v+k\X^p_{\gx v})(w+k\X^q_{\gx w}) - mr_{p,q}(k\X^p\otimes k\X^q_{\gex w})\\
&\subset& vw + vk\X^q_{\gx w} + k\X^p_{\gx v}k\X^q + k\X^q_{\gex w}k\X^p.\\
\end{eqnarray*}
Using the definition of the lexicographical order we see that the first and second subspace are contained $k\X^{l(u)}_{\gx u}$. For the third subspace take $a\in k\X^q_{\gex w},b\in k\X^p$. Then (because $u$ is Lyndon) $a\gex w\gx u$ and because $a$ is shorter than $u$ we obtain $ab\gx u$. Thus
\[[u]_r \in u + k\X^{p+q}_{\gx u}.\]
Now assume $u$ is not Lyndon. Let $u=u_1\ldots u_t$ be the Lyndon decomposition and let $v:=u_1, w:=u_2\ldots u_t$ and $p:=l(v),q:=l(w)$. Then
\begin{eqnarray*}
[u]_r &=& [v]_r[w]_r \in (v+k\X^p_{\gx v})(w+k\X^q_{\gx w})\\
&\subset& vw + v k\X^q_{\gx w} + k\X^p_{\gx v}k\X^q \subset u + k\X^{p+q}_{\gx u}.
\end{eqnarray*}
For the second part observe that if $r$ is diagonal then $[u]_r$ is just a linear combination of words $u'$ that are obtained from the word $u$ by permuting the letters of $u$. If also $v'$ is obtained form $v$ by permuting we have $\gamma_{u',v'} = \gamma_{u,v}$. Thus the diagonal part of $c$ acts on every monomial in $[u]_r\otimes [v]_r$ by multiplication with $\gamma_{u,v}$. Together with the preceeding lemma and the first part this completes the proof.
\qed\end{pf}

\section{The comultiplication in the tensor algebra}
\label{section_comult}
Now in preparation of the final theorem we will prove some combinatorical properties of the comultiplication of the tensor algebra of a left triangular braided vector space. So for this section fix a finite dimensional braided vector space $(V,c)$, assume that it is left triangular with respect to the basis $X$ and denote by $d$ the diagonal component of $c$. Abbreviate $[-] := [-]_{d^{-1}}$ using the inverse of the diagonal component (As $c$ is bijective it is easy to prove that the diagonal coefficients of $c$ are not zero).
\par The technical lemmas in this section are motivated by similar calculations done in \cite{kh.b} in the case of diagonal braidings. A key idea is to use the iterated commutators induced by the inverse of the diagonal component of the braiding. Just using $c$ or its inverse is not possible. The following lemma - a generalization of \cite[lemma 8]{kh.b} - requires this. As Kharchenko works with diagonal braiding he can use the inverse of the braiding itself for the commutator. It is a central observation that it does not matter in fact which admissible endomorphism one uses for the commutator. If the final theorem is proved for some admissible endomorphism $r$ it generalizes easily to any other admissible endomorphism (see remark \ref{remark_basis}).
\begin{lem}
\label{lemma_comultlyndon}
Let $u\in\X$ be a Lyndon word and $n := l(u)$. Then
\[ \Delta([u]) \in [u]\obar 1 + 1\obar [u] +  \lsum{{i+j=n}\atop{i,j\neq 0}}{}k\X_{\gx u}^i\obar k\X^j. \]
\end{lem}
\begin{pf}
Induction on $n=l(u)$. For $n=1$ nothing has to be proved. Assume $n>1$ and let $u=vw$ be the Shirshov decomposition of $u$. By induction we have
\[ \Delta([v]) \in [v]\obar 1 + 1\obar [v] + \lsum{{i+j=l(v)}\atop{i,j\neq 0}}{} k\X^i_{\gx v}\obar k\X^j \mtxt{and}\]
\[ \Delta([w]) \in [w]\obar 1 + 1\obar [w] + \lsum{{l+m=l(w)}\atop{l,m\neq 0}}{} k\X^l_{\gx w}\obar k\X^m.\]

Now we obtain
\begin{eqnarray*}
\Delta([v])\Delta([w])&\in&\left([v]\obar 1 + 1\obar [v] + \lsum{{i+j=l(v)}\atop{i,j\neq 0}}{} k\X^i_{\gx v}\obar k\X^j\right)\cdot\\
&&\left( [w]\obar 1 + 1\obar [w] + \lsum{{l+m=l(w)}\atop{l,m\neq 0}}{} k\X^l_{\gx w}\obar k\X^m\right)\\
&\subset& [v][w]\obar 1 + [v]\obar[w] + 1\obar [v][w] + \lsum{{i+j=n}\atop{i,j\neq 0}}{} k\X^i_{\gx u}\obar k\X^j
\end{eqnarray*}
using the following facts:\newline
$([v]\obar 1)(k\X^i_{\gx w}\obar k\X^j) \subset [v]k\X^i_{\gx w}\obar k\X^j\subset k\X^{i+l(v)}_{\gx u}\obar k\X^j$ by definition of the lexicographical order.
As $w\gx u$ as $u$ is Lyndon we have
\[(1\obar[v])([w]\obar 1) = c([v]\obar [w]) \in k\X^{l(w)}_{\gex w}\obar k\X^{l(v)}\subset k\X^{l(w)}_{\gx u}\obar k\X^{l(v)}.\]
Furthermore $(1\obar [v])(k\X^l_{\gx w}\obar k\X^m)\subset k\X^l_{\gx w}\obar k\X^{i+l(v)}\subset k\X^l_{\gx u}\obar k\X^{i+l(v)}$ using the same argument. Finally for all $i,j,l,m\in\N$ with $i<l(v)$ we have $(k\X^i_{\gx v}\obar k\X^j)(k\X^l\obar k\X^m) \subset k\X^i_{\gx v}k\X^l\obar k\X^{j+m}\subset k\X^{i+l}_{\gx u}\obar k\X^{j+m}$ because if $a\in\X$ is shorter than $v$ and $a>v$ then for all $b\in\X$ also $ab>vb$ .
\par
On the other hand
\begin{eqnarray*}
\Delta([w])\Delta([v]) &\in& \left([w]\obar 1 + 1\obar [w] + \lsum{{i+j=l(w)}\atop{i,j\neq 0}}{} k\X^i_{\gx w}\obar k\X^j\right)\cdot\\
&&\left( [v]\obar 1 + 1\obar [v] + \lsum{{l+m=l(v)}\atop{l,m\neq 0}}{} k\X^l_{\gx v}\obar k\X^m\right)\\
&\subset& [w][v]\obar 1 + d([w]\obar[v]) + 1\obar [w][v] + \lsum{{i+j=n}\atop{i,j\neq 0}}{} k\X^i_{\gx u}\obar k\X^j
\end{eqnarray*}
using $([w]\obar 1)(1 \obar [v]) = [w]\obar [v] \in k\X^{l(w)}_{\gex w}\obar k\X^{l(v)}\subset k\X^{l(w)}_{\gx u}\obar k\X^{l(v)}$ because again $w\gx u$ and $([w]\obar 1)(k\X^i_{\gx v}\obar k\X^j)\subset [w]k\X^i_{\gx v}\obar k\X^j \subset k\X^{l(w)+i}_{\gx wv}\obar k\X^j\subset k\X^{i+l(w)}_{\gx u}\obar k\X^j$ because $u$ is Lyndon and thus $wv\gx vw=u$.\newline
 Next $(1\obar [w])([v]\obar 1)=c([w]\obar [v])\in d([w]\obar [v]) + k\X^{l(w)}_{\gx w}\obar k\X^{l(v)}\subset d([w]\obar [v]) + k\X^{l(w)}_{\gx u}\obar k\X^{l(v)}$.
Finally $k\X_{\gx w}, k\X^i_{\gx v}\subset k\X_{\gx u}$ if $i<l(v)$ and $k\X_{\gx u}\obar k\X$ is a right ideal in $k\X\obar k\X$.
As $d([w]\obar [v]) = \gamma_{w,v} [v]\obar [w]$ we obtain
\begin{eqnarray*}
 \Delta([u]) &=& \Delta([v])\Delta([w])-\gamma_{w,v}{}^{-1}\Delta([w])\Delta([v])\\
&\in& 1\obar [u] + [u]\obar 1 + \lsum{{i+j=n}\atop{i,j\neq 0}}{} k\X^i_{\gx u}\obar k\X^j.
\end{eqnarray*}
\qed\end{pf}
To describe the comultiplication on arbitrary words we need an other subset of $\X$.
\begin{defn}
For $u,v\in\X$, $u$ a Lyndon word we write $v\mgt u$ if $u$ is smaller than the first Lyndon letter of $v$. Furthermore $\X_{\mgt u} := \{v\in\X | v\mgt u\}, \X_{\mgt u}^m := \X^m\cap\X_{\mgt u}$.
\end{defn}
Considering this subset is one of the key ideas in the step from the setting of diagonal braidings to that of triangular braidings. It can not be found in the work of Kharchenko. First we collect some auxiliary statements.
\begin{rem} Let $u,v\in\X$, $u$ a Lyndon word. Then:
\label{remark_mgt}
\begin{enumerate}
\item $\X_{\mgt u} = \{v\in\X | \mtxt{for all} i\in\N: v\gx u^i \}$.
\item If also $v$ is Lyndon $v\gx u$ implies $v\mgt u$.
\item If $v\mgt u$, then $\X_{\gex v}\subset\X_{\mgt u}$.
\item If $v$ is Lyndon, $v\gx u$ then $\X_{\gex v}\X\subset \X_{\mgt u}$.
\item If $v$ is Lyndon, $v\gx u$ then $\X_{\mgt v}\subset \X_{\mgt u}$.
\item If $v\gx u$ and $l(v)\leq l(u)$ then $v\mgt u$.
\item $\X_{\mgt u}\X \subset \X_{\mgt u}$.
\item If $i\in\N$ then $\X^{il(u)}_{\gex u^i} \X_{\mgt u} \subset \X_{\mgt u}$.
\end{enumerate}
\end{rem}
\begin{pf}
For part 1 assume first that $v\mgt u$ and let $v=v_1\ldots v_r$ be the Lyndon decomposition of $v$. This means that $v_1 \gx u$ by assumption and for all $i\in\N$ we obtain $v = v_1\ldots v_r \gx u^i$ by comparing the Lyndon letters lexicographically, keeping in mind that the Lyndon decomposition of $u^i$ consist of $i$ Lyndonletters $u$. On the other hand let $v\in \X$ with $v\gx u^i$ for all $i\in\N$. Again let $v=v_1\ldots v_r$ bet the Lyndon decomposition. Because of $v\gx u$ we have $v_1 \gex u$. Assume $v_1 = u$. Then we find $i\in\N$ such that $v_1=\ldots =v_i = u$ and $v_{i+1},\ldots,v_r \lx u$. If $i=r$ we have $v=u^r\lx u^{r+1}$, a contradiction. If $i<r$ we have $v=u^iv_{i+1}\ldots v_r \lx u^{i+1}$ by comparing the Lyndon letters - again a contradiction. Thus $v_1 \gx u$ and $v\mgt u$.\newline
Part 2 follows from the definition. Part 3: Let $w\in\X,w\gex v$. Then $w\gex v\gx u^i$ for all $i\in\N$ and so $w\in\X_{\mgt u}$. For part 4 consider $a,b\in\X$ with $a\gex v$. Then $ab \gex a \gex v \gx u^i$ for all $i\in\N$. Part 5 is trivial.\newline
Part 6: If $v$ is Lyndon, this is part 2. Otherwise let $v=v_1\ldots v_r$ be the Lyndon decomposition of $v$. Then $v_1\gex u$ and $l(v_1)<l(v)\leq l(u)$. So $v_1\gx u$ and $v\mgt u$. 
\newline
For part 7 let $a\in \X_{\mgt u},b\in\X$. Let $c$ be the first Lyndon letter of $a$. So we have $a\gex c\gx u$ and $ab\in \X_{\gex c}\X \subset \X_{\mgt u}$ by part 4. Part 8: Assume $a\gex u^i,l(a)=il(u)=l(u^i),b\mgt u$. First assume $j\geq i$. Then $ab \gex u^ib$ (because $a$ and $u^i$ have the same length) and $b\gx u^{j-i}$ and together $ab\gex u^ib\gx u^iu^{j-i}=u^j$. Now assume $j<i$. If $u^j$ is the beginning of $a$, then also of $ab$ and thus $ab\gx u^j$. If otherwise $u^j$ is not the beginning of $a$, then $u^j \lx a$ (because $u^i \lex a$) implies $u^j \lx ab$. In any case $ab \gx u^j$ and thus $ab\mgt u$.
\qed\end{pf}
\begin{rem}
Let $u\in\X$, u a Lyndon word and $p,q\in\N$. Then
\[ c(k\X^p\obar k\X_{\mgt u}^q) \subset k\X_{\mgt u}^q\obar k\X^p. \]
\end{rem}
\begin{pf}
Let $v,w\in\X,v\mgt u, l(v)=q,l(w)=p$. Then 
\[c(w\obar v) \in k\X_{\gex v}^q\obar k\X^p\subset k\X_{\mgt u}^q\obar k\X^p.\]
The last inclusion is by part 3 of the preceding remark.
\qed\end{pf}

\begin{cor}
Let $u\in\X$ be a Lyndon word, $n=l(u)$. Then
\[\Delta([u]) \in [u]\obar 1 +1\obar [u] +\lsum{{i+j=n}{}\atop{i,j\neq 0}}{} k\X^i_{\mgt u}\obar k\X^j. \]
\end{cor}
\begin{pf}
By lemma \ref{lemma_comultlyndon} because $\Delta$ is graded we obtain
\[\Delta([u]) \in [u]\obar 1 +1\obar [u] +\lsum{{i+j=n}\atop{i,j\neq 0}}{} k\X^i_{\gx u}\obar k\X^j, \]
but for $i<l(u)=n$ we have by part 6 of remark \ref{remark_mgt} $k\X^i_{\gx u}\subset k\X^i_{\mgt u}$. 
\qed\end{pf}

\par Thus up to terms of a special form (simple tensors whose left tensorand is made up of monomials having a first Lyndon letter bigger than $u$) the $[u]$ behave like primitive elements. The aim of this section is to extend this observation. In this spirit the next two lemmas are generalizations of calculations used in \cite{kh.b} to our situation. The more general context asks for a more careful formulation of statements and proofs.

\begin{lem}
\label{lemma_comultpotenz}
Let $v\in\X$ be a Lyndon word, $r\in\N$ and $n:=l(v^r)$. Then
\[\Delta([v^r]) \in \lsum{i=0}{r} \qgauss{r}{i}{\gamma_{v,v}} [v]^i\obar [v]^{r-i} + \lsum{{i+j=n}\atop{i,j\neq 0}}{} k\X^i_{\mgt v}\obar k\X^j.\]
\end{lem}

\begin{pf}
We use induction on $r$. The case $r=1$ is the preceeding corollary. So assume $r>1$. Then
\begin{eqnarray*}
\Delta([v^r]) &=& \Delta([v]^{r-1})\Delta([v])\\
&\in&\left(\lsum{i=0}{r-1} \qgauss{r-1}{i}{\gamma_{v,v}} [v]^i\obar [v]^{r-1-i} + \lsum{{i+j=n-l(v)}\atop{i,j\neq 0}}{} k\X^i_{\mgt v}\obar k\X^j\right)\\
&&\left([v]\obar 1 + 1\obar [v] + \lsum{{l+m=l(v)}\atop{l,m\neq 0}}{} k\X^l_{\mgt v}\obar k\X^m\right).
\end{eqnarray*}
Now note that $[v]^ik\X_{\mgt v}\subset k\X_{\mgt v}$ by part 8 of remark \ref{remark_mgt} and that thus $k\X_{\mgt v}\obar k\X$ is stable under left multiplication with elements from $[v]^i\obar k\X$. Furthermore $k\X_{\mgt v}\obar k\X$ is a right ideal in $k\X\obar k\X$.
As (again by part 8 of remark \ref{remark_mgt}) 
\begin{eqnarray*}
([v]^i\obar[v]^j)([v]\obar 1) &=& \gamma_{v,v}^j [v]^{i+1}\obar [v]^j + k\X^{l(v^{i+1)}}_{\gx v^{i+1}}\obar k\X^{l(v^k)}\\
&\subset& \gamma_{v,v}^j [v]^{i+1}\obar [v]^j + k\X^{l(v^{i+1)}}_{\mgt v}\obar k\X^{l(v^k)}
\end{eqnarray*}
this implies
\[\Delta([v^r])\in \lsum{i=0}{r-1}\qgauss{r-1}{i}{\gamma_{v,v}}\left(\gamma_{v,v}^i[v]^{i+1}\obar [v]^{r-1-i} + [v]^i\obar [v]^{r-i}\right) +  \lsum{{i+j=n}\atop{i,j\neq 0}}{} k\X^i_{\mgt v}\obar k\X^j\]
and using the recursion formula for the $q$-binomial coefficients we obtain the claim.
\qed\end{pf}

\begin{lem}
\label{lemma_comult}
Let $u\in\X$ and $u=u_1\ldots u_t v^r, u_1\gex\ldots\gex u_t\gx v$ be the Lyndon decomposition with $r,t\geq 1$. Define $z := u_1\ldots u_t, n:=l(u)$. Then
\[ \Delta([u]) \in [u]\obar 1 + \lsum{i=0}{r} \qgauss{r}{i}{\gamma_{v,v}} \gamma_{z,v}^i [v]^{i}\obar[z][v]^{r-i} + \lsum{{i+j=n}\atop{i,j\neq 0}}{} k\X^i_{\mgt v}\obar k\X^j.\]
\end{lem}
\begin{pf}
We use induction on $t$. First assume $t=1$. Then $z$ is Lyndon and $z\gx v$. By part 5 of remark \ref{remark_mgt} we have $k\X_{\mgt z}\subset k\X_{\mgt v}$. So we obtain using the preceeding lemma
\begin{eqnarray*}
\Delta([zv^r]) &=& \Delta([z])\Delta([v]^r)\\
&\in& \left([z]\obar 1 + 1\obar [z] + \lsum{{i+j=l(z)}\atop{i,j\neq 0}}{} k\X^i_{\mgt v}\obar k\X^j\right)\\
&&\left( \lsum{i=0}{r}\qgauss{r}{i}{\gamma_{v,v}} [v]^i\obar [v]^{r-i} + \lsum{{l+m=l(v^r)}\atop{l,m\neq 0}}{} k\X^l_{\mgt v}\obar k\X^m\right).
\end{eqnarray*}
Again $k\X_{\mgt v}\obar k\X$ is a right ideal and stable under left multiplication with $1\obar [z]$. \newline
Moreover by part 8 of remark \ref{remark_mgt} $(1\obar[z])([v]^i\obar[v]^{r-i})\in \gamma_{z,v}^i [v]^i\obar [z][v]^{r-i} + k\X^{l(v^i)}_{\gx v^i}\obar k\X^{n-l(v^i)}\subset  \gamma_{z,v}^i [v]^i\obar [z][v]^{r-i}+  k\X^{l(v^i)}_{\mgt v}\obar k\X^{n-l(v^i)}$. As furthermore $[z] \in k\X_{\gex z}\subset k\X_{\mgt v}$ (part 3 of remark \ref{remark_mgt}) we obtain
\[\Delta([u][v^r])\in [u][v]^r\obar 1 + \lsum{i=0}{r}\qgauss{r}{i}{\gamma_{v,v}} \gamma_{z,v}^i [v]^i\obar[z][v]^{r-i} + \lsum{{i+j=n}\atop{i,j\neq 0}}{} k\X^i_{\mgt v}\obar k\X^j.\]
So assume now $t>1$ and let $w := u_2\ldots u_t$. The induction hypothesis and the preceeding lemma give
\begin{eqnarray*}
\Delta([u_1]) &\in& [u_1]\obar 1 + 1\obar [u_1] + \lsum{{i+j=l(u)}\atop{i,j\neq 0}}{} k\X^i_{\mgt v}\obar k\X^j\\
\Delta([z][v]^r)&\in& [z][v^r]\obar 1 + \lsum{i=0}{r}\qgauss{r}{i}{\gamma_{v,v}} \gamma_{w,v}^i [v]^i\obar[w][v]^{r-i} \\
&&+ \lsum{{l+m=l(v^r)}\atop{l,m\neq 0}}{} k\X^l_{\mgt v}\obar k\X^m.
\end{eqnarray*}
Multiplying these we use that $k\X_{\mgt v}\obar k\X$ is a right ideal and stable under left multiplication with elements from $k\X_{\gx v}\obar k\X$ and $1\obar k\X$. Note $[u_1]\in k\X_{\mgt v}$, $[z][v^r]\in k\X_{\mgt v}$. Together with 
\[(1\obar [u_1])([v]^i\obar [z][v]^{r-i}) \in \gamma_{u_1,v}^i [v]^i\obar [u_1][z][v]^{r-i} + k\X^{l(v^i)}_{\mgt v}\obar k\X^{n-l(v^i)}\]
 we have
\begin{eqnarray*}
 \Delta([u_1][w][v^r])&\in& [u_1][w][v^r]\obar 1 + \lsum{i=0}{r}\qgauss{r}{i}{\gamma_{v,v}} \gamma_{w,v}^i\gamma_{u_1,v}^i [v]^i\obar [u_1][w][v]^{r-i}\\&& +\lsum{{i+j=n}\atop{i,j\neq 0}}{} k\X^i_{\mgt v}\obar k\X^j
\end{eqnarray*}
establishing the claim.
\qed\end{pf}

\section{A PBW basis for braided Hopf algebras with left triangular braiding}
\label{section_PBWbasis}
We will use the results of the preceding section to prove the existence of a PBW type basis made up of iterated skew commutators for a class of braided Hopf algebras. We start with a general definition of the term PBW basis.

\begin{defn}
\label{defi_pbwallg}
Let $A$ be an algebra, $P,S\subset A$ subsets and $h:S\rightarrow \N\cup\{\infty\}$. Let $(S,<)$ be a totally ordered set. Let $B(P,S,<,h)$ be the set of all products
\[ s_1^{e_1}\ldots s_t^{e_t}p\]
with $t\in\N_0$,$s_1>\ldots>s_t$, $s_i\in S$, $0<e_i<h(s_i)$ and $p\in P$. This set is called the \emph{PBW set} generated by $P$, $(S,<)$ and $h$. $h$ is called the height function of the PBW set.\newline
We say $(P,S,<,h)$ is \emph{PBW basis} of $A$ if $B(P,S,<,h)$ is a basis of $A$.
\end{defn}
Of course every algebra $A$ has the trivial PBW basis with $S=\emptyset$ and $P$ a basis of $A$. In the (braided) bialgebra case we are interested in the case that $P$ is a basis of the coradical. Thus in this section we are interested in the case where $P=\{1\}$. We will say that $B(S,<,h) := B(\{1\},S,<,h)$ is the PBW set (resp. PBW basis) generated by $(S,<)$ and $h$.\par

As in the preceding section fix a finite dimensional braided vector space $(V,c)$ that is left triangular with repect to a basis $X$ of $V$. Let $d$ be the diagonal component of $c$ and abbreviate $[-]:=[-]_{d^{-1}}$. Identify $T(V)$ with $k\X$.
\begin{defn}
Define the \textsl{standard order} in the following way. For two elements $u,v\in\X$ write $u\ges v$ if and only if $u$ is shorter than $v$ or if  $u\gex v$ lexicographically and $l(u)=l(v)$.
\end{defn}
In this order the empty word $1$ is the maximal element. As $X$ is assumed to be finite, this order fulfils the ascending chain condition, making way for inductive proofs. Define $\X_{\gs u}, \X_{\ges u}$ etc. in the obvious way. \newline
Now we will define the PBW set that will lead to the PBW basis of our braided Hopf algebra. The sets $S_I$ resp. $B_I$ are analogues of the sets of ``hard superletters'' resp. of ``monotonous words in hard superletters'' found in \cite{kh.b}. 
\begin{defn}
\label{defi_pbw} Let $I\subsetneq k\X$ be a biideal. Let $S_I$ be the set of Lyndon words from $\X$ that do not appear as (standard-) smallest monomial in elements of $I$:
\[ S_I := \{u\in\X |\mtxt{u is a Lyndon word and} u\notin k\X_{\gs u} + I \}.\]
For $u\in S_I$ define the height of $h_I(u)\in \{2,3,\ldots,\infty\}$ by
\[ h_I(u) := \operatorname{min}\{t\in\N | u^t\in k\X_{\gs u^t} +I\}\]
and let $B_I:=B(S_I,<,h_I)$ be the PBW set generated by $(S_I,<)$ and $(h_I)$, where $<$ denotes the lexicographical order.
\end{defn}

If $r$ is an admissible endomorphism of $V\otimes V$ and $U\subset \X$ is any subset define $[U]_r := \{[u]_r|u\in U\}$. Denote by $k[U]_r$ the $k$-linear subspace of $k\X$ spanned by $[U]_r$ (To avoid confusion with the notation for polynomial rings let me note that no polynomial rings will be considered during this section).\par 
This section will mainly be devoted to the proof of the following central theorem. Note that in the special case of diagonal braidings this theorem together with lemma \ref{lemma_heights} is a braided analogue of \cite[theorem 2]{kh.b}.
\begin{thm}
\label{thm_pbw}
Let $(V,c)$ be a finite dimensional braided vector space that is left triangular with respect to some basis $X$. Identify $T(V)$ with $k\X$ and let $I\subsetneq k\X$ be a braided biideal, $\pi:k\X\rightarrow (k\X)/I$ the quotient map. Then $\pi(B_I)$ and $\pi([B_I]_c)$ are bases of $(k\X)/I$.
\newline These are the truncated PBW bases generated by $\pi(S_I)$ resp. $\pi([S_I]_c)$ with heights $h_I(u)$ for $u\in S_I$.
\end{thm}

\begin{rem}
\label{rem_losekh}
The reader should observe that in changing from diagonal to triangular braidings we lost some information on the basis. Kharchenko shows that in the diagonal case every reducible Lyndon word $u$ is (modulo $I$) a linear combination of
\begin{itemize}
\item words of the same degree as $u$ that are nonascending products in PBW generators lexicographically smaller than $u$ and
\item words of degree smaller than that of $u$ that are nonascending products in arbitrary PBW generators.
\end{itemize}
It is an open question wether this (or something similar) can be done for triangular braidings.
\end{rem}

\subsection{Proof of theorem \ref{thm_pbw}}
We will omit the index corresponding to the subspace $I$ during the technical parts of this section. So we can introduce some new notation: for $n\in\N$ and a Lyndon word $v\in\X$ define $B^n := B\cap \X^n,B_{\mgt v} := B\cap\X_{\mgt v}$ and $B^n_{\mgt v} := B^n\cap B_{\mgt v}$.
\newline
The next proposition collects some statements which will be useful in the sequel. Analogues of Parts 1,3 and 4 are also used in \cite{kh.b}.
\begin{prop}
\label{prop_inclusions}
Let $r$ be an admissible endomorphism of $V\otimes V$. For every $m\in\N_0,u,v\in\X$, $v$ a Lyndon word we have the following inclusions
\begin{enumerate}
\item $k\X_{\ges u} \subset k[B_{\ges u}]_r + I$,
\item $k\X_{\mgt v}^m \subset k[B_{\mgt v}^m]_r + \lsum{0\leq i<m}{} k[B^i]_r + I$,
\item $k\X^m \subset \lsum{0\leq i\leq m}{} k[B^i]_r + I$,
\item $k\X = k[B]_r + I$.
\end{enumerate}
\end{prop}

\begin{pf}
First note that for $x,y,a,b\in\X$ $a\gs b$ implies $xay\gs xby$.
For part 1 we proceed by downward induction along the standard order (this works because the standard order satisfies the ascending chain condition). For $u=1$ the inclusion is valid. Now assume $u\ls 1$ and that for all words $\gs u$ the inclusion is valid. Let $w\ges u,m:=l(w)$. If $w\in B$ we have $w\in [w]_r + k\X_{\gs w}\subset k[B_{\ges w}]_r + I$ by induction. Assume $w\not\in B$ and let $w=w_1^{e_1}\ldots w_t^{e_t}$ be the Lyndon decomposition of $w$. As $w\not\in B$ we find an $0\leq i\leq t$ such that either $w_i\not\in S$ or $e_i\geq h(w_i)$. In the first case we have $w_i\in k\X_{\gs w_i}+I$, in the second case $w_i^{e_i}\in k\X_{\gs w_i^{e_i}}+I$. Anyway this implies $w\in k\X_{\gs w}+I$, but thus $[w]_r\in k\X_{\gs w}+I \subset k[B_{\gs w}]_r+I\subset k[B_{\gs u}]_r +I$ by induction.
Now for part 2 assume $w\in k\X^m_{\mgt v}$. Then by part 1 $w\in k\X_{\ges w}\subset k[B_{\ges w}]_r + I \subset k[B^m_{\gex w}]_r + \lsum{0\leq i<m}{}k[B^i]_r + I$. In view of \ref{remark_mgt} we have $k\X^m_{\gex w}\subset k\X^m_{\mgt v}$, finishing the proof.
For part 3 let $u_0$ be the smallest word of degree $m$. For $u\in \X^m$ we have $u\ges u_0$ and thus by part 1 $u\in k[B_{\ges u_0}]_r + I$ and this is a subset of $\lsum{0\leq i\leq m}{} k[B^i]_r$.
Finally Part 4 follows from Part 3.
\qed\end{pf}

For the rest of the proof we use the main ideas of \cite{kh.b}, but in a different and more general setting. As the triangular braiding requires a more careful analysis we work in the tensor algebra rather than in the quotient of the tensor algebra by the ideal $I$ (as Kharchenko does, not regarding the biproduct with the group algebra). This enables us to use linear maps as tools where Kharchenko argues by inspection of the occuring terms, a method for which our situation seems to be too complicated.

\begin{lem}
Let $r$ be an admissible endomorphism of $V\otimes V$. Then the set $[B_I]_r$ is linearily independent.
\end{lem}
\begin{pf}
The map $[-]_r:k\X\rightarrow k\X$ is homogenous. Furthermore it is surjective (use part 4 of the proposition above for the subspace $(0)$ to see that $k\X = k[\X]_r$). As the homogenous components are finite dimensional, $[-]_r$ is bijective and maps the linearily independent set $B_I$ onto $[B_I]_r$. So $[B_I]_r$ is linearily independent.
\qed\end{pf}

The next theorem is the key step to the final theorem combining the results of the preceeding section on the comultiplication with the results of this section.

\begin{thm}
\label{basisofcomplement}
Let $I\subsetneq k\X$ be a braided biideal in $k\X$. Then $[B_I]$ spans a $k$-linear complement of $I$.
\end{thm}
\begin{pf}
By the proposition above all we need to show is that $k[B_I]$ and $I$ have trivial intersection. For $n\geq 0$ let $U_n := \operatorname{k-span}\{[u]|u\in B_I, l(u)\leq n\}$. We show by induction on $n$ that for all $n\in\N$ we have $U_n \cap I = (0)$. First let $n=0$. Then $U_0 = k1$ and $I$ is proper ideal, so $U_0\cap I= (0)$. Now assume $n>0$. Assume $0\neq T\in U_n\cap I$. So we can write $T$ as a (finite) sum
\[ T = \lsum{{u\in B_I}\atop{l(u)\leq n}}{} \alpha_u [u].\]
We may assume that there is a $u\in B^n$ such that $\alpha_u\neq 0$. Now choose $v$ as the (lexicographically) smallest Lyndon letter occuring in the Lyndon decomposition of words $u\in B_I^n$ with $\alpha_u\neq 0$. Because of the minimality of $v$, it occurs in Lyndon decompositions of words $u\in B_I^n$ with $\alpha_u\neq 0$ only at the end. Let $t$ be the maximal number of occurences of $v$ in a Lyndon decomposition of word $u\in B_I^n$ with $\alpha_u\neq 0$. Thus we can decompose the sum for $T$ in the following way
\[ T = \lsum{u\in O}{}\alpha_u [a_u] [v]^t + \lsum{u\in P}{} \alpha_u [a_u] [v]^{t_u} + \lsum{u\in Q}{} \alpha_u [u] + \lsum{u\in R}{} \alpha_u [u],\]
where $O,P,Q,R\subset B_I$ and the words $a_u$ for $u\in O\cup P$ are chosen such that
\begin{itemize}
\item $O$ contains all words $u\in B_I^n$ of length $n$ with $\alpha_u\neq 0$ such that the Lyndon decomposition of $u$ ends with $v^t$. Furthermore $u=a_uv^t$.
\item $P$ contains all words $u\in B_I^n$ of length $n$ with $\alpha_u\neq 0$ such that the Lyndon decomposition of $u$ ends with $v^{t_u}$ for some $0\neq t_u<t$. Furthermore $u=a_uv^{t_u}$.
\item $Q$ contains all words $u\in B_I^n$ of length $n$ with $\alpha_u\neq 0$ that do not have the Lyndon letter $v$ in their Lyndon decomposition.
\item $R$ contains all words $u\in B_I^n$ of length less than $n$ with $\alpha_u\neq 0$.
\end{itemize}
Note that for all $u\in O$ we have $a_u\neq 1$ as $u\in B$. By analyzing the four terms we will show
\begin{eqnarray*}
 \Delta(T) &\in& T\obar 1 + [v]^t\obar\lsum{u\in O}{}\alpha_u \gamma_{a_u,v}^t [a_u] +\lsum{i=0}{t-1} [v]^i\obar k\X^{n-l(v^i)} +\\
&& \lsum{{i+j=n}\atop{i,j\neq n}}{}k[B^i_{\mgt v}]\obar k\X^j + \lsum{i+j<n}{} k[B^i]\obar k\X^j + I\obar k\X.
\end{eqnarray*}
At first consider $u\in O, l:= t$ or $u\in P,l:= t_u<t$. Note that then $a_u\neq 1$: $a_u=1$ would imply that $u=v^l\in k\X_{\gs v^l}$, but as well $l<h_I(v)$ because $u\in B$ which is a contradiction to the definition of $h_I(v)$. So we obtain
\begin{eqnarray*}
\Delta([a_u][v]^l) &\in& [a_u][v]^l\obar 1 + \lsum{i=0}{l}\qgauss{l}{i}{\gamma_{v,v}} \gamma_{a_u,v}^i [v]^i\obar [a_u][v]^{l-i} + \lsum{{i+j=n}\atop{i,j\neq 0}}{} k\X^i_{\mgt v}\obar k\X^j \\
&\subset& [a_u][v]^l\obar 1 + [v]^l\obar \gamma_{a_u,v}^l [a_u] + \lsum{i=0}{l-1} [v]^i\obar k\X^{n-l(v^i)} + \\
&& \lsum{{i+j=n}\atop{i,j\neq 0}}{} k[B^i_{\mgt v}]\obar k\X^j + \lsum{i+j<n}{} k[B^i]\obar k\X^j + I\obar k\X.
\end{eqnarray*}
In both cases ($u\in O$ or $u\in P$) this delivers the right terms in the sum for the word $u$. Now consider $u\in Q$ and let $w$ be the largest Lyndon letter occuring in the Lyndon decomposition of $u$. Then $u=a w^l$ for some $l\in\N,a\in\X$ and $w\gx v$ by construction of $v$. This leads to
\begin{eqnarray*}
\Delta([u]) &\in&[u]\obar 1 + \lsum{i=0}{l}\qgauss{l}{i}{\gamma_{w,w}}\gamma_{a,w}^i [w]^i\obar [a][w]^{l-i} + \lsum{{i+j=n}\atop{i,j\neq 0}}{} k\X^i_{\mgt w}\obar k\X^j\\
&\subset&[u]\obar 1+\lsum{{i+j=n}\atop{i,j\neq 0}}{} k\X^i_{\mgt v}\obar k\X^j\\
&\subset& [u]\obar 1 + \lsum{{i+j=n}\atop{i,j\neq 0}}{} k[B^i_{\mgt v}]\obar k\X^j + \lsum{i+j<n}{} k[B^i]\obar k\X^j + I\obar k\X.
\end{eqnarray*}
Finally consider $u\in R$. Then $l(u)<n$ and we obtain
\begin{eqnarray*}
\Delta([u]) &\in& [u]\obar 1 + \lsum{i+j<n}{} k\X^i\obar k\X^j\\
&\subset& [u]\obar 1 + \lsum{i+j<n}{}k[B^i]\obar k\X^j + I\obar k\X.
\end{eqnarray*}
Now by induction assumption we find a $\phi\in (k\X)^*$ such that
\[ \phi(I)=0, \phi([v]^t)=1\mtxt{and $\forall u\in B\setminus \{v^t\}$ with} l(u)<n: \phi([u])=0 . \]
With the inclusion showed above we get
\begin{eqnarray*}
(\phi\obar \id)\Delta(T)&\in& \lsum{u\in O}{} \alpha_u \gamma_{a_u,v}^t[a_u] + \lsum{{i+j=n}\atop{i,j\neq 0}}{} \phi(k[B^i_{\mgt v}])k\X^j +\\&& \lsum{i+j<n}{} \phi(k[B^i])k\X^j + \phi(I)k\X\\
&\subset& \lsum{u\in O}{} \alpha_u \gamma_{a_u,v}^t[a_u] + 0 + \lsum{j<n-tl(v)}{}k\X^j + 0\\
&\subset& \{k[B^{n-tl(v)}]\oplus U_{n-tl(v)-1})\setminus\{0\} \subset U_{n-1}\setminus \{0\}.
\end{eqnarray*}
Note that we can not obtain $0$, because we have a non-zero component in degree $n-tl(v)$.
On the other hand, as $I$ is a biideal, we have
\[ (\phi\obar \id)\Delta(T) \in \phi(I) k\X + \phi(k\X)I \subset I.\]
Thus $(\phi\obar \id)\Delta(T) \in I \cap (U_{n-1}\setminus\{0\})$, but by induction assumption this is the empty set, a contradiction.
\qed\end{pf}
\begin{cor}
\label{basisofcomplementbraided}
Let $I\subsetneq k\X$ be a braided biideal and $r$ an admissible endomorphism of $V\otimes V$. Then $[B_I]_r$ spans a $k$-linear complement of $I$.
\end{cor}
\begin{pf}
Again all we have to show is that $k[B]_r \cap I = \{0\}$. Assume $0\neq T\in k[B]_r\cap I$. We can write $T$ as
\[ T = \alpha [u]_r + \lsum{w\gs u}{}\beta_w [w]_r\]
with $\alpha\neq 0$. Then by the proposition \ref{prop_inclusions} and lemma \ref{lemma_smallestterm} we obtain first
\[ T \in \alpha u + k\X_{\gs u}\]
and from this
\[ T \in \alpha [u] + k\X_{\gs u} \subset \alpha [u] + k[B_{\gs u}] + I.\]
So now write $T=\alpha[u] + x + i$ with $x\in k[B_{\gs u}],i\in I$. Now by the theorem above we obtain $\alpha [u] + x\in I\cap k[B]= \{0\}$. Thus $\alpha = 0$ because $[B]_r$ is linearily independent, a contradiction.
\qed\end{pf}
Now theorem \ref{thm_pbw} follows as a special case of the following remark.
\begin{rem}
\label{remark_basis}
Let $I\subsetneq k\X$ be a braided biideal and $r$ an admissible endomorphism of $V\otimes V$. Then $\pi([B_I]_r)$ is a basis of $(k\X)/I$.
\end{rem}
\begin{pf}
By corollary \ref{basisofcomplementbraided} $[B_I]_r$ is a basis for a complement of $I$. As $I=\ker\pi$, $\pi$ induces a $k$-linear isomorphism $k[B_I]_r\rightarrow k\X/I$, mapping the basis $[B_I]_r$ into $\pi([B_I]_r)$. This proves the remark. The theorem follows by using $r=0$ and $r=c$.
\qed\end{pf}

\subsection{A result on the height function}
The next lemma is useful to ensure that all elements of $S_I$ have infinite height if $\operatorname{char}(k)=0$ and the diagonal coefficients of the braiding are powers of one element that is not a root of unity. Its analogue in the case of diagonal braidings is already contained in the definition of the height function in \cite{kh.b}.

\begin{lem}
\label{lemma_heights}
Let $(V,c)$ be a finite dimensional braided vector space that is left triangular with respect to some basis $X$. Identify $T(V)$ with $k\X$ and let $I\subsetneq k\X$ be a braided biideal and $v\in S_I$.
Define the scalar $\gamma_{v,v}\in k$ by
\[ d(v\otimes v) = \gamma_{v,v} v\otimes v\]
where $d$ is the diagonal component of $c$.\\[0,2cm]
If $h:=h_I(v)<\infty$, then $\gamma_{v,v}$ is a root of unity. In this case let $t$ be the order of $\gamma_{v,v}$. If $\operatorname{char}k = 0$ then $\gamma_{v,v}\neq 1$ and $h=t$. If $\operatorname{char}k = p > 0$ then $h=tp^l$ for some $l\in\N$.
\end{lem}

\begin{pf}
Let $n:=l(v^h)$. We have an element of $I$ of the form
\[ T:=[v]^h + \lsum{{u\gx v^h}\atop {l(u)=n}}{}\alpha_u [u] + \lsum{l(w)<n}{} \alpha_w [w]\in I.\]
 For every $u\in\X^n$ with $u\gx v^h$ we have $u\mgt v$ and thus we obtain for the coproduct using lemmas \ref{lemma_comultpotenz} and \ref{lemma_comult}
 \begin{eqnarray*}
 \Delta(T)&\in& T\obar 1 + \lsum{i=0}{h-1}\qgauss{h}{i}{\gamma_{v,v}} [v]^i\obar [v]^{h-i} + \lsum{{i+j=n}\atop{j\neq 0}}{} k\X^i_{\mgt v}\obar k\X^j + \lsum{{i+j<n}\atop {j\neq 0}}{} k\X^i\obar k\X^j\\
 &\subset&T\obar 1 + \lsum{i=0}{h-1}\qgauss{h}{i}{\gamma_{v,v}} [v]^i\obar [v]^{h-i} + \lsum{{i+j=n}}{} k[B^i_{\mgt v}]\obar k[B^j] + \\&& \lsum{{i+j<n}}{} k[B^i]\obar k[B^j] + I\obar k\X + k\X\obar I.
 \end{eqnarray*}

 Now because of $k\X\obar k\X = (I\obar k\X + k\X\obar I) \oplus (k[B]\obar k[B])$ we can construct a $k$-linear map $\phi_1:k\X\otimes k\X\rightarrow k$ such that
 \begin{eqnarray*}
&& \forall (b,b')\in (B\times B)\setminus\{(v,v^{h-1})\}: \phi_1([b]\obar [b']) = 0 ,\\
 &&\phi_1([v]\obar [v]^{h-1})=1,\\
 &&\phi_1(I\obar k\X + k\X\obar I) = 0.
 \end{eqnarray*}
 As $T\in I$ we have $\phi_1\Delta(T)=0$ and on the other hand using what we proved above
 \[0=\phi_1\Delta(T) = \qgauss{h}{h-1}{\gamma_{v,v}} = 1 + \gamma_{v,v} +\ldots+\gamma_{v,v}{}^{h-1}.\]
 This shows that $\gamma_{v,v}$ is a root of unity, say of order $t$ (set $t=1$ if $\gamma_{v,v}=1$). Let $p:= \operatorname{char} k$ and define $q$ by
 \[ q := \left\{\begin{array}{ccc}p &\mbox{if}&p>0\\1 & \mbox{if}&p=0.\end{array}\right.\]
Now we can write $h=tq^l a$ with $a,l\in\N$. If $q\neq 1$ we may assume that $q$ does not divide $a$. We want to show that $a=1$. So assume $a> 1$. In this case we can construct a $k$-linear map $\phi_2:k\X\otimes k\X\rightarrow k$ with
 \begin{eqnarray*}
&& \forall (b,b')\in (B\times B)\setminus\{(v^{tq^l},v^{tq^l(a-1)})\}: \phi_2([b]\obar [b']) = 0 ,\\
&& \phi_2([v]^{tq^l}\obar [v]^{tq^l(a-1)})=1,\\
&& \phi_2(I\obar k\X + k\X\obar I) = 0.
 \end{eqnarray*}
Using that $\gamma_{v,v}$ is a primitive $t$-th root of unity (resp. $\gamma_{v,v}=1$ and $t=1$) we obtain that in $k$
 \[ 0=\phi_2\Delta(T) = \qgauss{tq^la}{tq^l}{\gamma_{v,v}} = \qgauss{q^la}{q^l}{} = \qgauss{a}{1}{} = a.\]
This is a contradiction to the assumptions we made on $a$.
Thus $h=tq^l$. In particular if $\operatorname{char} k=0$, then $q=1$ and because $t=h>1$ we obtain $\gamma_{v,v}\neq 1$.
 \qed\end{pf}

\newpage

\section{Right triangular braidings}
\label{sect_right}
In principle one could do a similar proof as above for right triangular braidings, but an easy argument shows that the right triangular case follows from the left triangular case. Obviously $c$ is a right triangular braiding if and only if $\tau c\tau$ is a left triangular braiding, where $\tau$ denotes the usual flip map. The key observation is

\begin{prop}
\label{prop_opposite}
Let $(R,\mu,\eta,\Delta,\eps,c)$ be a braided bialgebra. \newline Then also $R^{op,cop} := (R,\mu\tau,\eta,\tau\Delta,\eps,\tau c \tau)$ is a braided bialgebra.
\end{prop}
\begin{pf}
Of course $(R,\mu^{op},\eta)$ is an algebra, $(R,\Delta^{cop},\eps)$ is a coalgebra and $(R,\tau c\tau)$ is a braided vector space. Checking the compatibility of $\mu^{op},\eta,\Delta^{cop},\eps$ with $\tau c\tau$ is tedious. We will do one example, namely the calculation that $\tau c\tau\circ (\mu^{op}\otimes R) = (R\otimes \mu^{op})(\tau c\tau\otimes R)(R\otimes\tau c\tau)$. We calculate
\begin{eqnarray*}
&&(R\otimes\mu^{op})(\tau c\tau\otimes R)(R\otimes\tau c\tau) = \\
&& (R\otimes\mu)\underline{(R\otimes\tau)(\tau\otimes R)(R\otimes\tau)}(R\otimes\tau)(c\tau\otimes R)(R\otimes\tau c\tau) = \\
&& \underline{(R\otimes\mu)(\tau\otimes R)(R\otimes\tau)}\circ\underline{(\tau\otimes R)(R\otimes\tau)(c\otimes R)}(\tau\otimes R)(R\otimes\tau c\tau) = \\
&& \tau(\mu\otimes R)\circ(R\otimes c)\underline{(\tau\otimes R)(R\otimes\tau)(\tau\otimes R)}(R\otimes \tau c\tau)=\\
&& \tau(\mu\otimes R)(R\otimes c)\underline{(R\otimes\tau)(\tau\otimes R)(R\otimes c)}(R\otimes\tau)=\\
&& \tau\underline{(\mu\otimes R)(R\otimes c)(c\otimes R)}\circ\underline{(R\otimes\tau)(\tau\otimes R)(R\otimes\tau)}=\\
&& \tau c \underline{(R\otimes\mu)\circ(\tau\otimes R)(R\otimes\tau)}(\tau\otimes R)=\\
&& \tau c\tau(\mu\otimes R)(\tau\otimes R) =\\
&& \tau c\tau(\mu^{op}\otimes R),
\end{eqnarray*}
where we use (in this order): $\tau^2=\id_{V\otimes V}$, the braid equation for $\tau$, $\mu,c$ commute with $\tau$, again the braid equation for $\tau$ and $\tau^2=\id_{V\otimes V}$, $c$ commutes with $\tau$, $\mu$ commutes with $c$ and the braid equation for $\tau$ and finally again that $\mu$ commutes with $\tau$. The other calculations work similarily (use graphical calculus as a tool for intiution).\newline
Finally we have to check that $\Delta:R\rightarrow R\obar R$ and $\eps:R\rightarrow k$ are algebra morphisms, where $R\obar R$ is an algebra with multiplication $(\mu\otimes\mu)(R\otimes \tau c\tau\otimes R)$. For $\eps$ this is clear. 
For $\Delta $ we have to check
\[ \Delta^{cop}\mu^{op} = (\mu^{op}\otimes\mu^{op})(R\otimes\tau c\tau\otimes R)(\Delta^{cop}\otimes\Delta^{cop}).\]
As $R$ is a braided bialgebra the left hand side is
\[ \tau\Delta\mu\tau = \tau(\mu\otimes\mu)(R\otimes c\otimes R)(\Delta\otimes\Delta)\tau.\]
Now because $\Delta,\mu$ commute with $\tau$ this is equal to
\[ (\mu\otimes\mu)(R\otimes \tau\otimes R)(\tau\otimes\tau)(R\otimes\tau c\tau\otimes R)(\tau\otimes\tau)(R\otimes \tau\otimes R)(\Delta\otimes\Delta).\]
Thus it suffices to show
\[ (R\otimes\tau\otimes R)(\tau\otimes\tau)(R\otimes\tau c\tau\otimes R)(\tau\otimes\tau)(R\otimes\tau\otimes R) = 
(\tau\otimes\tau)(R\otimes\tau c\tau\otimes R)(\tau\otimes\tau),\]
but this is trivial (check on elements).
\qed\end{pf}
Related material can be found in \cite{AG1}.

Assume now that $(V,c)$ is a braided vector space. Denote the braided tensor bialgebra defined in section \ref{section_brha} by $(T(V,c),\mu,\eta,\Delta_c,\eps,c)$. As an algebra this is $T(V)$.
\begin{prop}
\label{prop_isoopcop}
Let $(V,c)$ be a braided vector space. Let 
\[\phi:T(V,c)\rightarrow T(V,\tau c\tau)^{op,cop}\]
be the unique algebra morphism $T(V)\rightarrow T(V)^{op}$ given by $\phi|V=\id_V$. Then $\phi$ is an isomorphism of braided bialgebras. \par
For $v_1,\ldots,v_n\in V$ we have
\[\phi(v_1\otimes\ldots\otimes v_n)= v_n\otimes\ldots\otimes v_1.\]
\end{prop}
\begin{pf}
Lemma \ref{epiforbraidedHA} gives us the existence of $\phi$ as a morphism of braided bialgebras because $(V,c)$ is a braided subspace of $T(V,\tau c\tau)^{op,cop}$. By construction we see that $\phi$ is bijective and has the form given in the lemma.
\qed\end{pf}
Now we can prove the existence of the PBW basis in the right triangular case.
\begin{thm}
Assume $(V,c)$ is a finite dimensional right triangular braided vectorspace and $I\subsetneq T(V,c)$ is a braided biideal. Then there is a totally ordered subset $S\subset T(V,c)$ and a height function $h:S\rightarrow\N\cup\{\infty\}$ such that the images of the PBW set generated by $S$ and $h$ form a basis of $T(V,c)/I$.\newline
Let
\[ \phi:T(V,c)\rightarrow T(V,\tau c\tau)^{op,cop}\]
be the isomorphism from proposition \ref{prop_isoopcop}.
We have
\[ S = \phi^{-1}(S_{\phi(I)})\:\:\mbox{and}\:\:h = h_{\phi(I)}\phi\]
and the order on the set $S$ is the opposite of the order on $S_{\phi(I)}$.
\end{thm}
\begin{pf}
As $\phi(I)$ is a braided biideal in $T(V,\tau c\tau)^{op,cop}$ it is also a braided biideal in $T(V,\tau c\tau)$. As $c$ is right triangular we have that $\tau c \tau$ is left triangular. So we find a set $S_{\phi(I)}\subset T(V,\tau c\tau)$ with a total ordering $<$ and a height function $h_{\phi(I)}:S\rightarrow \N\cup\{\infty\}$ such that the PBW set generated by these data in $T(V,\tau c \tau)$ is a basis for a complement of $\phi(I)$. The PBW set generated in $T(V,\tau c\tau)^{op,cop}$ by $S_{\phi(I)}$ with reversed order and height function $h_{\phi(I)}$ is the same set and thus also a basis for a complement of $\phi(I)$. The claim follows by transfering this set to $T(V,c)$ via $\phi^{-1}$.
\qed\end{pf}

Together with the characterization of Nichols algebras given in \cite{AS5} we obtain
\begin{cor}
Let $(V,c)$ be a braided vectorspace. Then
\[\Nichols(V,\tau c \tau)^{op,cop}\isom \Nichols(V,c)\]
as braided graded Hopf algebras.
\end{cor}
\begin{pf}
In \cite{AS5} the kernel of the canonical projection
\[\pi_c:T(V,c)\rightarrow \Nichols(V,c)\]
is characterized as the sum over all coideals that are also homogenous ideals generated by elements of degree $\geq 2$. Using this we see that the map $\phi$ from proposition \ref{prop_isoopcop} maps $\ker\pi_c$ bijectively into $I:=\ker\pi_{\tau c\tau}\subset T(V,\tau c\tau)$. As obviously
\[ T(V,\tau c \tau)^{op,cop}/I = (T(V,\tau c\tau)/I)^{op,cop} \isom \Nichols(V,\tau c\tau)^{op,cop}\]
we get that $\phi$ induces the desired isomorphism of braided Hopf algebras.
\qed\end{pf}

\section{Application to pointed Hopf algebras with abelian coradical}
\label{sect_applHA}
In this section we will show how to obtain a PBW basis for a Hopf algebra generated by an abelian group and a finite dimensional $G$-module spanned by skew primitive elements. On one hand this yields a generalization of the result in \cite{kh.b} as there the skew primitive elements are assumed to be semi-invariants (i.e. that the group acts on them by a character). On the other hand we lose some properties of the basis as already mentioned in remark \ref{rem_losekh}.

Let $A=\cup_{n\geq 0} A_n$ be a filtered algebra. We can define a map
\[ \pi:A \rightarrow \gr A\]
by setting $\pi(0):= 0$ and for all $0\neq a\in A: \pi(a) := a+A_{n-1}$ for the unique $n\geq 0$ such that $a\in A_n\setminus A_{n-1}$ (where $A_{-1}:=\{0\}$ as usual). We will use this map to obtain PBW bases for $A$ from homogenous PBW bases of the associated graded algebra $\gr A$.

\begin{prop}
Let $A=\cup_{n\geq 0} A_n$ be a filtered algebra and $(P,S,<,h)$ a PBW basis for $\gr A$ such that $P\subset \gr A(0)=A_0$ and $S$ is made up of homogenous elements. Then there is a PBW basis $(P,S',<',h')$ of $A$ such that for all $a,b\in S'$
\[ \pi(a)\in S, h'(a) = h(\pi(a))\:\:\mbox{and}\:\:a<b\Leftrightarrow \pi(a)<'\pi(b).\]
\end{prop}
\begin{pf}
For all $s\in S\cap\gr A(n)$ we find $\hat{s}\in A_n\setminus A_{n-1}$ such that $\pi(\hat{s})=s$. Define
\[ S' := \{\hat{s}|s\in S\}.\]
The map $S\rightarrow S', s\mapsto \hat{s}$ is bijective. So we can transfer the height funtion $h$ and the order $<$ to $S'$ obtaining $h'$ and $<'$.\\[0,2cm]
Assume we have $b := s_1^{e_1}\ldots s_r^{e_r}p\in B(P,S,<,h)$. We define a lift
\[\hat{b} := \hat{s_1}^{e_1}\ldots\hat{s_r}^{e_r}p\in B(P,S',<',h'). \]
As the $s_i$ and $p$ are homogenous (say of degrees $n_i$ and $0$), also $b$ is homogenous, say of degree $n$. Then
\[ b = (\hat{s_1}+A_{n_1-1})^{e_1}\ldots(\hat{s_r}+A_{n_r-1})^{e_r}(p+A_{-1}) 
=\hat{s_1}^{e_1}\ldots\hat{s_r}^{e_r}p + A_{n-1} = \hat{b} + A_{n-1}\]
in $\gr A(n)=A_n/A_{n-1}$. We have $\hat{b}\in A_n\setminus A_{n-1}$, because otherwise ($\hat{b}\in A_{n-1}$) we had $b=0$, but this is an element of a basis.\\[0,2cm]
Let $B_n := B(P,S,<,h)\cap\gr A(n)$ and $\hat{B}_n := \{ \hat{b}|b\in B_n\}$. We will show by induction on $n\geq 0$ that $\hat{B}_0\cup\ldots\cup\hat{B}_n$ generates $A_n$ as a vector space. The case $n=0$ is trivial as $A_0=\gr A(0)$ has basis $P$. Assume $n\geq 0$ and $a\in A_n\setminus A_{n-1}$. We have $\pi(a) = a+A_{n-1}\in\gr A(n)$ and thus $\pi(a)$ is a linear combination of elements of $B_n$ i.e. 
\[\pi(a) = a + A_{n-1} \in kB_n = \sum\limits_{b\in B_n} k(\hat{b}  + A_{n-1})\]
So we get that $a$ is a linear combination of elements from $\hat{B}_n$ and $A_{n-1}$ and by induction assumption $a$ is a linear combination of elements from $\hat{B}_0\cup\ldots\cup\hat{B}_n$.\\[0,2cm]
We are left to show that $B(P,S',<',h')$ is linearily independend. Assume we have for all $b\in B_n$ scalars $\alpha_b\in k$ such that
\[ \sum\limits_{b\in B_n}\alpha_b \hat{b} \in A_{n-1}.\]
It suffices to show for all $b$ that $\alpha_b = 0$. As we have seen above we have for all $b\in B(P,S,<,h): b = \pi(\hat{b})$. Thus we have in $\gr A(n)$:
\[\sum\limits_{b\in B_n} \alpha_b b = 
\sum\limits_{b\in B_n} \alpha_b (\hat{b} + A_{n-1}) = 
\left(\sum\limits_{b\in B_n} \alpha_b \hat{b}\right) + A_{n-1} = 0\]
As $B_n$ is a linearily independend we obtain for all $b\in B_n\setminus B_{n-1}: \alpha_b=0$.
\qed\end{pf}

\begin{thm}
Let $k$ be an algebraically closed field.
Let $H$ be a Hopf algebra generated by an abelian group $G$ and skew primitive elements $a_1,\ldots,a_t$ such that the subvectorspace of $H$ spanned by $a_1,\ldots,a_t$ is stable under the adjoint action of $G$. Then $H$ has a PBW basis $(G,S,<,h)$.
\end{thm}

\begin{pf}
First we may assume that for all $1\leq i\leq t$
\[ \Delta(a_i) = g_i\otimes a_i + a_i\otimes 1.\]
Let $H_n$ be the subspace of $H$ generated by all products of elements of $G$ and at most $n$ factors from $\{a_1,\ldots,a_t\}$. This defines a Hopf algebra filtration of $H$. It is well known from \cite{AS5,Radford_HAproj,Majid_crossprod} that we can decompose the associated graded Hopf algebra
\[ \gr H \isom R\# kG,\]
as graded Hopf algebras, where $R$ is a braided graded Hopf algebra in ${}^{kG}_{kG}\mYD$ generated by the finite dimensional Yetter-Drinfel'd submodule $R(1)\subset P(R)$. As a $G$ module $R(1)$ is isomorphic to $ka_1+\ldots+ka_t$ with the adjoint $G$-action.



Example \ref{exmp_jordan} shows that the braiding on $R(1)$ is triangular because the group $G$ is abelian. So by the PBW theorem \ref{thm_pbw} we find a PBW basis $(\{1\},S,<,h)$ of $R$. This implies that $(1\#G,S\#1,<,h)$ is a PBW basis of $\gr H$ and thus we find a PBW basis of $H$ using the proposition above.
\qed\end{pf}

\section{Examples}
\label{section_examples}
Finally we will deal with some interesting examples mentioned in \cite{A1} that are not of diagonal type by remark \ref{remark_notdiag}, namely the Nichols algebras of simple \Uq{sl_2} modules of low dimension (and type $+1$). Let $q\in\C$ be not a root of unity and $(M,c) = L(n,+1)$ be the simple \Uq{sl_2} module of dimension $n+1$ and type $+1$. Denote its natural basis (see e.g. \cite{Jantzen}) by $x_0,\ldots,x_n$ and order this basis by $x_0<\ldots <x_n$. Then the braidings from example \ref{expl_uqg} are left (and right) triangular with respect to this basis. For calculations we assume for the map $f$ that $f(\frac{\alpha}{2},\frac{\alpha}{2})=q^{-2}=f(0,0)$ (see  \cite[chapters 3 and 7]{Jantzen}).
\begin{enumerate}
\item $n=1$: As for example shown in \cite{AS5}, $\mathcal{B}(M,c)$ is a quadratic algebra. The relation in degree two is
\begin{eqnarray*}  &x_0x_1 - qx_1x_0 = 0.&\end{eqnarray*}
Thus the set of PBW generators $S$ contains only Lyndon words in $x_0,x_1$ that do not have $x_0x_1$ as a subword. This means $S=\{x_0,x_1\}$. As the diagonal coefficients of $c$ are powers of $q$ by lemma \ref{lemma_heights} all elements have infinite height. Thus the elements of the form $x_1^ix_0^j,i,j\in\N_o$ form a basis for $\mathcal{B}(M,c)$.\par
\item $n=2$: Here we have the following relations in degree two:
\begin{eqnarray*}
& x_0x_1-q^2x_1x_0 = 0,&\\ 
&x_1x_2-q^2x_2x_1 = 0,&\\
&x_0x_2+(q^2-1)x_1x_1-x_2x_0=0. &
\end{eqnarray*}
So all words in $S$ are Lyndon words in $x_0,x_1,x_2$ that do not contain one of $x_0x_1,x_1x_2,x_0x_2$ as a subword. As surely $x_0,x_1,x_2\in S$ this leaves $S=\{x_0,x_1,x_2\}$. Again by lemma \ref{lemma_heights} all elements of $S$ have infinite height and thus the elements of the form $x_2^ix_1^jx_0^k,i,j,k\in\N_0$ form a basis of $\mathcal{B}(M,c)$. In particular this is a quadratic algebra.\par
\item $n=3$: In this case the space of relations of degree two can be calculated using Maple. It is generated by the relations
\begin{eqnarray*}
&x_0x_1-q^3x_1x_0=0,&\\
&x_0x_2+\frac{1-q^4}{q}x_1^2-x_2x_0=0,&\\
&q^3x_0x_3+q^2(q^2+1-q^4)x_1x_2+(q-q^3-q^5)x_2x_1-x_3x_0=0,&\\
&x_1x_3+\frac{1-q^6}{q(q^2+1)}x_2^2-x_3x_1=0,&\\
&x_2x_3-q^3x_3x_2=0.&
\end{eqnarray*}
By combining these relations one obtains the additional relations
\begin{eqnarray*}
&(q^4-q^2+1)x_1x_2x_2-q(q^6+1)x_2x_1x_2+(q^4-q^2+1)q^4x_2x_2x_1 = 0,&\\
&x_1x_1x_2-q(q^2+1)x_1x_2x_1+q^4x_2x_1x_1 = 0.&
\end{eqnarray*}
As $q$ is a root of unity, the leading coefficients in these relations are not zero: the zeros of $X^4-X^2+1$ are primitive $12$-th roots of unity as 
\[X^{12}-1 = (X^4-X^2+1)(X^2+1)(1-X^6).\]
Thus $S$ can only contain Lyndon words in $x_0,x_1,x_2,x_3$ that do not contain a subword from the following list:
\[ x_0x_1, x_0x_2, x_0x_3, x_1x_3, x_2x_3, x_1^2x_2, x_1x_2^2.\]
It follows that $S\subset \{x_0,x_1,x_2,x_3,x_1x_2\}$. None of these words can be expressed by standard-bigger ones as we can see from the relations of degree 2. Thus $S=\{x_0,x_1,x_2,x_3,x_1x_2\}$ and the elements of the form $x_3^ax_2^b(x_1x_2)^cx_1^dx_0^e, a,b,c,d,e\in\N_0$ form a basis of $\mathcal{B}(M,c)$. In particular $\mathcal{B}(M,c)$ is a quadratic algebra.

\end{enumerate}

\nocite{Humphreys}
\bibliography{../../tex/promotion.bib}

\begin{thebibliography}{10}
\expandafter\ifx\csname url\endcsname\relax
  \def\url#1{\texttt{#1}}\fi
\expandafter\ifx\csname urlprefix\endcsname\relax\def\urlprefix{URL }\fi

\bibitem{kh.b}
V.~Kharchenko, {A quantum Analog of the Poincar\'e-Birkhoff-Witt Theorem},
  Algebra and Logic Vol. 38~(4) (1999) 259--276.

\bibitem{Reutenauer}
C.~Reutenauer, {Free Lie Algebras}, Vol.~7 of London Mathematical Society
  Monographs, New Series, Clarendon Press, London, 1993.

\bibitem{LR}
P.~Lalonde, A.~Ram, {Standard Lyndon bases of Lie algebras and enveloping
  algebras}, Trans. Amer. Math. Soc. 347~(5) (1995) 1821--1830.

\bibitem{Lusztig_fdha}
G.~Lusztig, {Finite-dimensional Hopf algebras arising from quantized enveloping
  algebras}, J. Amer. Math. Soc. 3~(1) (1990) 257--296.

\bibitem{Lusztig}
G.~Lusztig, {Introduction to quantum groups}, no. 110 in Progress in
  Mathematics, Birkh\"auser, 1993.

\bibitem{Rosso_pbw}
M.~Rosso, {An analogue of the PBW theorem and the universal $R$-matrix for
  $U_hsl(N+1)$}, Comm. Math. Phys. 124~(2) (1989) 307--318.

\bibitem{Rosso_invent}
M.~Rosso, {Quantum Groups and quantum Shuffles}, Inventiones Math. 133 (1998)
  399--416.

\bibitem{Ringel}
C.~Ringel, {PBW-bases of quantum groups}, J. reine angewandte Mathematik 470
  (1990) 51--88.

\bibitem{Drin85}
V.~Drinfel'd, {Hopf algebras and the quantum Yang-Baxter equation}, Dokl. Akad.
  Nauk SSSR 283:5 (1985) 1060--1064.

\bibitem{Drin87}
V.~Drinfel'd, {Quantum groups}, in: Proc. Int. Cong. Math. (Berkeley 1986),
  Amer. Math. Soc., Providence, RI, 1987, pp. 798--820.

\bibitem{Jimbo85}
M.~Jimbo, {A $q$-differences analogue of $U(\mathfrak{g})$ and the Yang-Baxter
  equation}, Lett. Math. Phys. 10 (1985) 63--69.

\bibitem{Leclerc}
B.~Leclerc, {Dual canonical Bases, quantum Shuffles and $q$-Characters},
  preprint at arXiv:math.QA/0209133v1, to appear in Math. Zeitschrift .

\bibitem{A1}
N.~Andruskiewitsch, {Some Remarks on Nichols Algebras}, preprint at
  arXiv:math.QA/0301064 .

\bibitem{Lothaire}
M.~Lothaire, {Combinatorics on Words}, Vol.~17 of Encyclopedia of Mathematics
  and its Applications, Addison-Wesley, 1983.

\bibitem{Ufnarovski}
V.~Ufnarovski, {Combinatorial and asymptotic Methods in Algebra}, in:
  Kostrikin, Shafarevich (Eds.), Algebra VI, no.~57 in Encyclopaedia of
  Mathematical Sciences, Springer, 1994.

\bibitem{Takeuchi_survey}
M.~Takeuchi, {Survey of braided Hopf Algebras}, in: New trends in Hopf Algebra
  Theory, no. 267 in Contemporary Mathematics, 1999, pp. 301--323.

\bibitem{Sbg_borel}
P.~Schauenburg, {A Characterization of the Borel-like Subalgebras of quantum
  enveloping Algebras}, Comm. in Algebra 24 (1996) 2811--2823.

\bibitem{AS5}
N.~Andruskiewitsch, H.-J. Schneider, {Pointed Hopf Algebras}, in: New
  directions in Hopf algebras, Vol.~43 of MSRI Publications, Cambridge
  University Press, 2002, pp. 1--68.

\bibitem{Ich_pointed}
S.~Ufer, {Braidings and pointed Hopf algebras}, preprint 2003.

\bibitem{Jantzen}
J.~Jantzen, {Lectures on quantum Groups}, Vol.~6 of Graduate Studies in
  Mathematics, Am. Math. Soc., 1995.

\bibitem{kh.c}
V.~Kharchenko, {A combinatorical Approach to the Quantification of Lie
  Algebras}, Pacific J. Math. 203~(1) (2002) 191--233.

\bibitem{AG1}
N.~Andruskiewitsch, M.~Gra\~na, {Braided Hopf algebras over non-abelian
  Groups}, Bol. Acad. Ciencieas (Cordoba) 63 (1999) 45--78.

\bibitem{Radford_HAproj}
{D.E. Radford}, {Hopf algebras with projection}, J. Algebra 92 (1985) 322--347.

\bibitem{Majid_crossprod}
S.~Majid, {Crossed products by braided groups and bosonization}, J. Algebra 163
  (1994) 165--190.

\bibitem{Humphreys}
J.~Humphreys, {Reflection Groups and Coxeter Groups}, Vol.~29 of Cambridge
  studies in advanced mathematics, Cambridge University Press, 1994.

\end{thebibliography}
\end{document}